\documentclass[3p]{elsarticle}

\usepackage{amsthm}
\usepackage{amssymb}
\usepackage{amsmath}
\usepackage{amsfonts}
\usepackage{epsfig}
\usepackage[mathscr]{eucal}
\usepackage[ruled,vlined]{algorithm2e}
\usepackage{algorithmic}
\usepackage{setspace,color}
\usepackage{tikz}
\usepackage{hyperref}
\usepackage[capitalize]{cleveref}
\usepackage{numcompress}
\usepackage[caption=false]{subfig}
\usepackage{multirow}

\usepackage[numbers]{natbib}
\usepackage{stackrel}

\usetikzlibrary{snakes}
\usetikzlibrary{arrows,shapes}

\newtheorem{thm}{Theorem}[section]
\newtheorem{proposition}[thm]{Proposition}
\newtheorem{lemma}[thm]{Lemma}
\newtheorem{theorem}[thm]{Theorem}

\newdefinition{definition}[thm]{Definition}
\newdefinition{rmk}[thm]{Remark}
\newdefinition{notation}[thm]{Notation}
\newdefinition{example}[thm]{Example}

\newproof{pf}{Proof}
\newproof{poth1}{Proof of \cref{Th1}}
\newproof{poth3}{Proof of \cref{Th3}}

\numberwithin{equation}{section}

\def\R{{\mathbb R}}

\def\N{{\mathbb N}}

\def\Z{{\mathbb Z}}

\def\PP{\mathbb{P}}
\def\EE{\mathbb{E}}
\def\BB{{\mathbb B}}

\def\su{\mathrm{supp}}

\def\Om{\Omega}

\def\ep{\epsilon}
\def\s{\sigma}
\def\f{\frac}

\def\na{\nabla}
\def\la{\langle}
\def\ra{\rangle}

\def\a{{\boldsymbol a}}

\def\rd{{\mathrm d}}
\def\e{{\boldsymbol e}}

\def\bsf{{\boldsymbol f}}

\def\bsg{{\boldsymbol g}}
\def\bsh{{\boldsymbol h}}
\def\bi{{\mathbf i}}

\def\bk{{\boldsymbol k}}

\def\bn{\boldsymbol{n}}
\def\bp{\boldsymbol{p}}
\def\bq{\boldsymbol{q}}

\def\bu{\boldsymbol{u}}

\def\bsr{{\boldsymbol r}}

\def\x{\boldsymbol{x}}

\def\z{{\boldsymbol z}}

\def\bS{{\mathbf S}}

\def\mA{{\mathcal A}}

\def\mE{{\mathcal E}}

\def\mI{{\mathcal I}}

\def\mN{{\mathcal N}}

\def\mP{{\mathcal P}}

\def\mR{{\mathcal R}}

\def\mW{{\mathcal W}}

\def\bmA{{\boldsymbol \mA}}

\def\bmI{{\boldsymbol \mI}}

\def\bmP{{\boldsymbol \mP}}

\def\bmR{{\boldsymbol \mR}}

\def\bmW{{\boldsymbol \mW}}

\def\msE{{\mathscr E}}
\def\msF{{\mathscr F}}

\def\msM{{\mathscr M}}

\def\msR{{\mathscr R}}

\def\msX{{\mathscr X}}

\def\bmsF{{\boldsymbol \msF}}

\def\aal{\boldsymbol{\alpha}}
\def\bbe{\boldsymbol{\beta}}
\def\dde{\boldsymbol{\delta}}

\def\lam{\boldsymbol{\lambda}}
\def\oom{\boldsymbol{\omega}}

\def\nnu{\boldsymbol{\nu}}

\def\bi{\begin{itemize}} \def\ei{\end{itemize}}
\def\be{\begin{eqnarray*}}
\def\ee{\end{eqnarray*}}

\def\0{{\mathbf 0}}

\newcommand{\beq}{\begin{equation}}
\newcommand{\eeq}{\end{equation}}

\def\xxi{{\boldsymbol{\xi}}}
\def\eet{{\boldsymbol\eta}}
\def\zze{{\boldsymbol\zeta}}

\def\La{{\boldsymbol\Lambda}}

\def\Ps{{\boldsymbol \Psi}}

\def\wt{\widetilde}
\def\wh{\widehat}

\def\Na{\boldsymbol \nabla}

\def\bOm{\boldsymbol{\Omega}}

\def\la{\langle}
\def\ra{\rangle}

\def\XXint#1#2#3{{\setbox0=\hbox{$#1{#2#3}{\int}$ }
\vcenter{\hbox{$#2#3$ }}\kern-.55\wd0}}

\crefformat{equation}{\textup{#2(#1)#3}}
\crefrangeformat{equation}{\textup{#3(#1)#4--#5(#2)#6}}
\crefmultiformat{equation}{\textup{#2(#1)#3}}{ and \textup{#2(#1)#3}}
{, \textup{#2(#1)#3}}{, and \textup{#2(#1)#3}}
\crefrangemultiformat{equation}{\textup{#3(#1)#4--#5(#2)#6}}%
{ and \textup{#3(#1)#4--#5(#2)#6}}{, \textup{#3(#1)#4--#5(#2)#6}}{, and \textup{#3(#1)#4--#5(#2)#6}}

\Crefformat{equation}{#2Equation~\textup{(#1)}#3}
\Crefrangeformat{equation}{Equations~\textup{#3(#1)#4--#5(#2)#6}}
\Crefmultiformat{equation}{Equations~\textup{#2(#1)#3}}{ and \textup{#2(#1)#3}}
{, \textup{#2(#1)#3}}{, and \textup{#2(#1)#3}}
\Crefrangemultiformat{equation}{Equations~\textup{#3(#1)#4--#5(#2)#6}}%
{ and \textup{#3(#1)#4--#5(#2)#6}}{, \textup{#3(#1)#4--#5(#2)#6}}{, and \textup{#3(#1)#4--#5(#2)#6}}

\crefdefaultlabelformat{#2\textup{#1}#3}

\title{Approximation Theory of Wavelet Frame Based Image Restoration}

\author[HKUST,HKUSTSZ]{Jian-Feng Cai\fnref{JFC}}
\ead{jfcai@ust.hk}
\author[TJU]{Jae Kyu Choi\fnref{JKC}\corref{cor}}
\ead{jaycjk@tongji.edu.cn}
\author[HHU]{Jianbin Yang\fnref{JBY}}
\ead{jbyang@hhu.edu.cn}

\cortext[cor]{Corresponding author}

\address[HKUST]{Department of Mathematics, The Hong Kong University of Science and Technology, Hong Kong}
\address[HKUSTSZ]{HKUST Shenzhen-Hong Kong Collaborative Innovation Research Institute, Futian, Shenzhen, China.}
\address[TJU]{School of Mathematical Sciences, Tongji University, Shanghai, 200092 China}
\address[HHU]{Department of Mathematics, Hohai University, Nanjing, 211100 China}

\fntext[JFC]{J. F. Cai is supported in part by Hong Kong Research Grants Council (HKRGC) GRF grants 16309518, 16309219, 16310620, and 16306821, and by the Project of Hetao Shenzhen-Hong Kong Science and Technology Innovation Cooperation Zone (HZQB-KCZYB-2020083).}
\fntext[JKC]{J. K. Choi is supported in part by the National Natural Science Foundation of China (NSFC) Youth Program 11901436, the Shanghai Science and Technology Committee 20JC1413500, and the Fundamental Research Funds for the Central Universities.}
\fntext[JBY]{J. Yang is supported in part by the National Natural Science Foundation of China (NSFC) 11771120 and the Fundamental Research Funds for the Central Universities B200202005.}

\begin{document}

\begin{abstract} In this paper, we analyze the error estimate of a wavelet frame based image restoration method from degraded and incomplete measurements. We present the error between the underlying original discrete image and the approximate solution which has the minimal $\ell_1$-norm of the canonical wavelet frame coefficients among all possible solutions. Then we further connect the error estimate for the discrete model to the approximation to the underlying function from which the underlying image comes.
\end{abstract}


\maketitle


\pagestyle{myheadings}
\thispagestyle{plain}
\markboth{Jian-Feng Cai, Jae Kyu Choi, and Jianbin Yang}{Approximation Theory of Wavelet Frame Based Image Restoration}

\section{Introduction}\label{Introduction}

In many practical image restoration tasks, observed images or measurements are often incomplete in the sense that features of interest in an image are missing partially and/or corrupted by noise. The corresponding inverse problem is to restore the underlying image $\bsf$ from its incomplete measurements given by
\begin{align}\label{Linear_IP_Incomplete}
\bsg[\bk]=\left\{\begin{array}{cl}
(\bmA\bsf)[\bk]+\eet[\bk],~&\bk\in\La,\vspace{0.5em}\\
\text{unknown},~&\bk\in\bOm\setminus\La.
\end{array}\right.
\end{align}
In \cref{Linear_IP_Incomplete}, $\eet$ is a measurement error which is assumed to be uncorrelated with $\left(\bmA\bsf\right)\big|_{\La}$, and $\bmA$ is some linear operator describing the process of imaging, acquisition, and communication. As convention, we consider discrete images as real-valued functions on a two dimensional regular grid $\bOm\subseteq\Z^2$, and $\La$ is a subset of $\bOm$ where $\bsg$ is available or reliable. Nevertheless, the discussions on other dimensions will be almost the same with slight modifications. Note that this paper involves both functions and their discrete counterparts. We shall use regular characters to denote functions and use bold-faced characters to denote their discrete analogies. For example, we use $u$ as an element in a function space, while we use $\bu$ to denote its corresponding discretized version (the type of discretization will be made clear later).

The measurements could be part of (degraded) images, transmitted data, and sensor data. The restoration of missing data from a given incomplete measurement is an essential part of imaging science area where the final image is utilized for the visual interpretation or for the automatic analysis \cite{Y.R.Li2011}. In the literature, numerous methods have been proposed under many different settings, e.g. \cite{M.Bertalmio2003,J.F.Cai2010,T.F.Chan2002,M.Elad2005} for image inpainting, \cite{T.F.Chan2006a,Y.W.Wen2012,X.Ye2013,X.Zhang2010a} for wavelet domain inpainting, \cite{T.Goldstein2009,W.Guo2014,M.Lustig2007} for compressed sensing MRI restoration, \cite{J.Choi2016,B.Dong2013a} for sparse view/angle CT restoration, and \cite{E.J.Candes2006,R.H.Chan2005} for miscellaneous applications. That being said, we forgo to give detailed surveys on these various applications and the interested reader should consult the references mentioned above for details. Instead, assuming that
\begin{align}\label{NoiseAssumption}
\mE\left(\eet\right)=\f{1}{|\La|}\sum_{\bk\in\La}\eet[\bk]=0~~~\text{and}~~~\mathrm{Var}\left(\eet\right)=\f{1}{|\La|}\sum_{\bk\in\La}\eet[\bk]^2=\eta^2,
\end{align}
we focus on establishing the approximation analysis of the following wavelet frame analysis based approach \cite{J.F.Cai2009/10,M.Elad2005,J.L.Starck2005}:
\begin{align}\label{WaveletFrameModel}
\min_{\bu\in[0,M]^{\bOm}}\left\|\lam\cdot\bmW\bu\right\|_1~~~\text{subject to}~~\f{1}{|\La|}\sum_{\bk\in\La}\left|(\bmA\bu)[\bk]-\bsg[\bk]\right|^2\leq\eta^2.
\end{align}
In \cref{WaveletFrameModel}, $\bmW$ is some wavelet frame transform, and $\left\|\lam\cdot\bmW\bu\right\|_1$ is the weighted $\ell_1$-norm of the wavelet frame coefficients of $\bu$ which promotes the regularity. Since the key features of image are in general modelled as singularities of piecewise smooth function, we only consider the images of low regularity throughout this paper. For the clarity of presentation, the detailed definition and explanation on this weighted $\ell_1$-norm will be postponed until \cref{ErrorAnalysis}.

The settings of \cref{Linear_IP_Incomplete} considered in this paper are as follows:
\begin{enumerate}
\item $\bOm=\left\{0,1,\cdots,N-1\right\}^2$, where $N\in\N$,
\item $\La\subsetneq\bOm$, $|\La|=m$, and $\La$ is uniformly randomly chosen from all $m$-subsets of $\bOm$,
\item $\s_{\min}(\bmA)>0$ and both $\s_{\min}(\bmA)$ and $\left\|\bmA\right\|_{\infty}$ are independent of $|\bOm|$, where $\s_{\min}(\bmA)$ denotes the smallest singular value of $\bmA$ and $\left\|\bmA\right\|_{\infty}$ denotes the operator norm induced by the $\ell_{\infty}$ norm.
\end{enumerate}
We denote by $\rho:=m/|\bOm|$ the density of the pixels available. In addition, we assume that the operator $\bmA$, the measurement $\bsg$ and the error $\eet$ are given and fixed, even though $\bmA$ and $\eet$ can be viewed as particular realizations of some random variables, e.g. the Gaussian or Bernoulli sensing matrix for $\bmA$ and the i.i.d. white Gaussian noise for $\eet$. Hence, the only random variable we consider is the data set $\La$ which is uniformly randomly drawn from all $m$-subsets of $\bOm$. Related examples include the wavelet inpainting \cite{T.F.Chan2006a,Y.W.Wen2012,X.Ye2013,X.Zhang2010a} to restore an image from a degraded and randomly missed data due to the unreliable transmission channel, the simultaneous deblurring and inpainting due to e.g. simultaneous out-of-focusing (or camera shaking) and dust spots or cracks in film \cite{H.Ji2011}, and the image restoration in the presence of Gaussian and impulsive noise \cite{L.Bar2007,L.Bar2006,J.F.Cai2008,J.F.Cai2010a,J.F.Cai2009,R.H.Chan2005,H.Ji2011,Y.R.Li2011}, etc.

Obviously, \cref{WaveletFrameModel} takes the wavelet frame based image inpainting ($\bmA=\bmI$ being the identity operator) whose approximation property is studied in \cite{J.F.Cai2011} as a special case. Thus, our goal is to extend the approximation analysis in \cite{J.F.Cai2011} to the generic image restoration tasks \cref{Linear_IP_Incomplete}. To begin with, we assume that $\bsf\in[0,M]^{\bOm}$ and $\left\|\lam\cdot\bmW\bsf\right\|_1<\infty$. The first condition enforces the range of pixel values on $\bsf$, which is typically set to be $M=255$ for the grayscale images, and the second condition enforces some regularity on $\bsf$. Let $\bu^{\La}$ be the solution to \cref{WaveletFrameModel}, which is possible as \cref{WaveletFrameModel} admits at least one solution with the minimal weighted $\ell_1$-norm of wavelet frame coefficients subject to the constraint \cite{G.Aubert2006}. For $\rho=1$, we have
\begin{align*}
\f{1}{|\bOm|}\left\|\bu^{\La}-\bsf\right\|_{\ell_2(\bOm)}^2\leq\f{2\s_{\min}^{-2}(\bmA)}{|\bOm|}\left(\left\|\bmA\bu^{\La}-\bsg\right\|_{\ell_2(\bOm)}^2+\left\|\bsg-\bmA\bsf\right\|_{\ell_2(\bOm)}^2\right)\leq 4\s_{\min}^{-2}(\bmA)\eta^2,
\end{align*}
which means that if $\rho=1$ and $\eta=0$, the unique solution to \cref{WaveletFrameModel} is the original image $\bsf$. Hence, we aim to focus on analyzing what happens when $\rho<1$ and $\eta>0$. More precisely, we will show that, under some mild assumptions, the following inequality
\begin{align}\label{Goal}
\f{1}{|\bOm|}\left\|\bu^{\La}-\bsf\right\|_{\ell_2(\bOm)}^2\leq C_1\rho^{-1/2}|\bOm|^{-\gamma}\left(\log_2|\bOm|\right)^{3/2}+C_2\eta^2
\end{align}
holds with probability at least $1-|\bOm|^{-1}$. In \cref{Goal}, $\gamma$ is a positive constant related to the regularity of $\bsf$, and $C_1$ and $C_2$ are constants independent of $|\bOm|$, $\rho$, and $\eta$. Briefly speaking, as long as the data set is sufficiently large, one has a pretty good chance to restore $\bsf$ by solving \cref{WaveletFrameModel}.

In the literature, there are efficient numerical algorithms developed to solve \cref{WaveletFrameModel} with a guaranteed convergence to minimizer \cite{J.F.Cai2009/10,B.Dong2012,H.Ji2011}. These numerical algorithms aim to find a sparse approximation of the underlying image using a properly chosen wavelet frame transform. Then a sparsity promoting regularization term (such as the widely used $\ell_1$-norm) is used to regularizing smooth image components while preserving image singularities such as edges and ridges. Hence, one may consider the approximation analysis of the numerical algorithms in the context of compressed sensing (CS) (e.g. \cite{E.J.Candes2006,E.J.Candes2005,E.J.Candes2006a,Donoho2006}). However, we would like to mention that our setting of \cref{WaveletFrameModel} is different from that of CS. First of all, we only assume the decay of wavelet frame coefficients $\bmW\bsf$ in the sense that $\left\|\lam\cdot\bmW\bsf\right\|_1$ is bounded rather than the explicit sparsity of $\bsf$ and/or $\bmW\bsf$, so that we can impose a low regularity assumption on $\bsf$. Second and most importantly, the ``sensing matrix'' of \cref{Linear_IP_Incomplete} does not satisfy the conditions required for the theoretical analysis of CS. In our setting, the sensing matrix is $\bmR_{\La}\bmA$, where $\bmR_{\La}$ denotes the projection onto $\La$; $\bmR_{\La}\bsf[\bk]=\bsf[\bk]$ for $\bk\in\La$, and $\bmR_{\La}\bsf[\bk]=0$ for $\bk\in\bOm\setminus\La$. Since $\La$ is randomly chosen from all $m$-subsets of $\bOm$ with a given fixed $m$, $\bmR_{\La}$ does not satisfy the concentration inequality in \cite{H.Rauhut2008}, and $\bmR_{\La}\bmA$ may not satisfy the restricted isometry property condition in general, all of which mean that there may not be information in the measurement $\bsg$ sufficient for the exact restoration \cite{J.F.Cai2011}. Instead, the main tools are the uniform law of large numbers and an estimation of covering number \cite{F.Cucker2002} of the hypothesis space, which are widely used in classical empirical process \cite{A.W.vanderVaart1996} and statistical learning theory \cite{Vapnik1998}. We further mention that our error analysis uses the newly involved estimation for covering number. More precisely, even though we also use the special structure of the set and the max-flow min-cut theorem in graph theory to estimate covering number, we improve the estimate in \cite[Theorem 2.4]{J.F.Cai2011} by relaxing the constraint on the radius in the covering number.


As a byproduct of the approximation analysis, we can further connect the approximation property of \cref{WaveletFrameModel} to the asymptotic approximation to the underlying function in terms of the shift invariant space generated by a scaled compactly supported basis function. Such an asymptotic approximation is extensively studied in terms of the Strang-Fix condition \cite{I.Daubechies2003,B.Dong2021,M.J.Johnson2009,J.Yang2017}. Notice that these asymptotic approximations require a certain order of regularity of underlying function in terms of e.g. the Sobolev space \cite{Adams1975}, even in the case of $\ell_1$ minimization approach in \cite{B.Dong2021,J.Yang2017}. Meanwhile, our asymptotic approximation uses the decay condition of wavelet frame coefficients of underlying function and the Bessel property of the shifts of scaled basis function only, as an extension of \cite{J.F.Cai2011} for the wavelet frame based image inpainting to the generic image restoration tasks.

The rest of this paper is organized as follows. The explicit formulation of \cref{Goal} is presented in \cref{ErrorAnalysis}. More precisely, we begin with introducing basic preliminaries of wavelet frames in \cref{Preliminaries}. The error estimates for \cref{WaveletFrameModel} in the discrete setting is given in \cref{WaveletFrameAnalysis}, and in \cref{FunctionApproximation}, we further establish the connection to function approximations in terms of the multiresolution approximation. Finally, all technical details are presented in \cref{TechnicalProofs}, and \cref{Conclusion} concludes this paper with relevant future works.

\section{Approximation error of wavelet frame based image restoration}\label{ErrorAnalysis}

\subsection{Preliminaries of wavelet frames}\label{Preliminaries}

We begin with a brief introduction on wavelet tight frames. For details, one may consult \cite{Daubechies1992,I.Daubechies2003,A.Ron1997} for theories of frames and wavelet frames, \cite{Shen2010} for a short survey on the theory and applications of frames, and \cite{B.Dong2013,B.Dong2015} for more detailed surveys.

A countable set $\left\{\varphi_n:n\in\Z\right\}\subseteq L_2(\R^d)$ with $d\in\N$ is called a tight frame of $L_2(\R^d)$ if
\begin{align}\label{TightFrame}
\left\|f\right\|_{L_2(\R^d)}^2=\sum_{n\in\Z}\left|\left\la f,\varphi_n\right\ra\right|^2~~~~~\text{for all}~~~f\in L_2(\R^d),
\end{align}
where $\la\cdot,\cdot\ra$ is the inner product on $L_2(\R^d)$, and $\la f,\varphi_n\ra$ is called the canonical coefficient of $f$.

For given $\Ps=\left\{\psi_1,\cdots,\psi_r\right\}\subseteq L_2(\R^d)$ and $J\in\N$, the corresponding quasi-affine system $\msX^J(\Ps)$ generated by $\Ps$ is defined by the collection of the dilations and the shifts of the elements in $\Ps$:
\begin{align}\label{QASystem}
\msX^J(\Ps)=\left\{\psi_{\alpha,l,\bk}:1\leq\alpha\leq r,~l\in\Z,~\bk\in\Z^d\right\}
\end{align}
with $\psi_{\alpha,l,\bk}$ being defined as
\begin{align}\label{QAFramelet}
\psi_{\alpha,l,\bk}=\left\{\begin{array}{cl}
2^{\f{ld}{2}}\psi_{\alpha}(2^l\cdot-\bk)~~&~~l\geq J;\vspace{0.5em}\\
2^{\left(l-\f{J}{2}\right)d}\psi_{\alpha}(2^l\cdot-2^{l-J}\bk)~~&~~l<J.
\end{array}\right.
\end{align}
When $\msX^J(\Ps)$ forms a tight frame of $L_2(\R^d)$, each $\psi_1,\cdots,\psi_r$ is called a (tight) framelet and the entire system $\msX^J(\Ps)$ is called a (tight) wavelet frame. In particular when $J=0$, we simply write $\msX(\Ps)=\msX^0(\Ps)$. Note that in the literature, the affine system is widely used, which corresponds to the decimate wavelet (frame) transform. The quasi-affine system, which corresponds to the undecimated wavelet (frame) transformation, was first introduced and analyzed in \cite{A.Ron1997}. Throughout this paper, we only discuss the quasi-affine system \cref{QAFramelet} because it generally performs better in image restoration than the widely used affine system, and the connection to the underlying function in continuum is more natural \cite{J.F.Cai2012,J.F.Cai2016,B.Dong2017a}. The interested reader can find further details on the affine wavelet frame systems and its connections to the quasi-affine wavelet frame systems in \cite{A.Chai2007,B.Dong2013,A.Ron1997}.

The constructions of framelets $\Ps$, which are desirably (anti-)symmetric and compactly supported functions, are usually based on a multiresolution analysis (MRA) generated by some refinable function $\phi$ with a refinement mask $\a_0$ such that
\begin{align}\label{MRA-RF}
\phi=2^d\sum_{\bk\in\Z^d}\a_0[\bk]\phi(2\cdot-\bk).
\end{align}
The idea of an MRA based construction of $\Ps=\left\{\psi_1,\cdots,\psi_r\right\}\subseteq L_2(\R^d)$ is to find finitely supported masks $\a_l$ such that
\begin{align}\label{MRA-Fra}
\psi_{\alpha}=2^d\sum_{\bk\in\Z^d}\a_{\alpha}[\bk]\phi(2\cdot-\bk)~~~~~\alpha=1,\cdots,r.
\end{align}
The sequences $\a_1,\cdots,\a_r$ are called wavelet frame mask or the high pass filters of the system, and the refinement mask $\a_0$ is also called the low-pass filter.

The unitary extension principle (UEP) of \cite{A.Ron1997} provides a general theory of the construction of MRA based tight wavelet frames. Briefly speaking, as long as $\left\{\a_0,\a_1,\cdots,\a_r\right\}$ are compactly supported and their Fourier series
\begin{align*}
\wh{\a}_{\alpha}(\xxi)=\sum_{\bk\in\Z^d}\a_{\alpha}[\bk]e^{-i\xxi\cdot\bk},~~~~~~\alpha=0,\cdots,r,~~~\xxi\in\R^d
\end{align*}
satisfy
\begin{align}\label{UEP}
\sum_{\alpha=0}^r\left|\wh{\a}_{\alpha}(\xxi)\right|^2=1~~~~\text{and}~~~~\sum_{\alpha=0}^r\wh{\a}_{\alpha}(\xxi)\overline{\wh{\a}}_{\alpha}(\xxi+\nnu)=0
\end{align}
for all $\nnu\in\left\{0,\pi\right\}^d\setminus\left\{\0\right\}$ and $\xxi\in[-\pi,\pi]^d$, the quasi-affine system $\msX(\Ps)$ with $\Ps=\left\{\psi_1,\cdots,\psi_r\right\}$ defined by \cref{MRA-Fra} forms a tight frame of $L_2(\R^d)$, and the filters $\left\{\a_0,\a_1,\cdots,\a_r\right\}$ form a discrete tight frame on $\ell_2(\Z^d)$ \cite{B.Dong2013}.

One of the most widely used examples is the piecewise linear B-spline \cite{I.Daubechies2003} for $L_2(\R)$, which has one refinable function and two framelets with the associated filters
\begin{align*}
\a_0=\f{1}{4}\big[\begin{array}{ccc}
1&2&1
\end{array}\big],~~~~\a_1=\f{\sqrt{2}}{4}\big[\begin{array}{ccc}
1&0&-1
\end{array}\big],~~\text{and}~~\a_2=\f{1}{4}\big[\begin{array}{ccc}
-1&2&-1
\end{array}\big].
\end{align*}
Indeed, it can be shown that the above $\left\{\a_0,\a_1,\a_2\right\}$ satisfies \cref{UEP}, so that $\msX(\Ps)$ with $\Ps=\left\{\psi_1,\psi_2\right\}$ defined by \cref{MRA-Fra} forms a tight frame on $L_2(\R)$.

One possible way to construct a tight frame on $L_2(\R^2)$ is by taking tensor products of univariate tight framelets \cite{J.F.Cai2012,J.F.Cai2016,Daubechies1992,B.Dong2013}. Given a set of univariate masks $\big\{\a_0,\a_1,\cdots,\a_r\big\}$, we define two dimensional masks $\a_{\aal}[\bk]$ with $\aal=(\alpha_1,\alpha_2)$ and $\bk=(k_1,k_2)$ as
\begin{align*}
\a_{\aal}[\bk]=\a_{\alpha_1}[k_1]\a_{\alpha_2}[k_2],~~~~~~0\leq\alpha_1,\alpha_2\leq r,~~\bk=(k_1,k_2)\in\Z^2
\end{align*}
so that the corresponding two dimensional refinable function and framelets are defined as
\begin{align*}
\psi_{\aal}(\x)=\psi_{\alpha_1}(x_1)\psi_{\alpha_2}(x_2),~~~~~0\leq\alpha_1,\alpha_2\leq r,~~\x=(x_1,x_2)\in\R^2
\end{align*}
with $\psi_0=\phi$ for convenience. If the univariate masks $\left\{\a_l:l=0,\cdots,r\right\}$ are constructed from UEP, then it can be verified that $\left\{\a_{\aal}:\aal\in\{0,\cdots,r\}^2\right\}$ satisfies \cref{UEP} and thus $\msX(\Ps)$ with
\begin{align*}
\Ps=\left\{\psi_{\aal}:\aal\in\{0,\cdots,r\}^2\setminus\{\0\}\right\}
\end{align*}
forms a tight frame for $L_2(\R^2)$.

The advantage of framelet is that we can easily derive the discrete tight frame system by framelet decomposition and reconstruction algorithms of \cite{I.Daubechies2003}. Throughout this paper, we only consider the MRA based tensor product wavelet frame system on $\ell_2(\bOm)$ with $\bOm=\left\{0,\cdots,N-1\right\}^2$ for simplicity. For $l\geq0$ we construct filters $\a_{l,\aal}$ by
\begin{align}\label{FilterDefinition}
\a_{l,\aal}=\wt{\a}_{l,\aal}\ast\wt{\a}_{l-1,\0}\ast\cdots\ast\wt{\a}_{0,\0}~~\text{with}~~\wt{\a}_{l,\aal}[\bk]=\left\{\begin{array}{cl}
\a_{\aal}[2^{-l}\bk],&\bk\in 2^l\Z^2;\vspace{0.5em}\\
0,&\bk\notin 2^l\Z^2,
\end{array}\right.
\end{align}
where $\a_{0,\aal}=\a_{\aal}$, and $\ast$ denotes the standard discrete convolution on $\Z^2$.

Let $\bu\in\ell_2(\bOm)$. To define an appropriate discrete framelet transform, or an analysis operator (see, e.g., \cite{B.Dong2013}), we will impose the periodic boundary conditions. Throughout this paper, we identify $\ell_p(\bOm)$ with the space of $\ell_p$ sequences on $\Z^2$ with fundamental period on each variable to be $N$. We introduce the following periodization operator $\bmP_J:\ell_1(\Z^2)\to\ell_1(\bOm)$:
\begin{align*}
\bmP_N\left(\bu\right)[\bk]=\sum_{\bk'\in\Z^2}\bu[\bk+N\bk']~~~~~~~\bk\in\bOm.
\end{align*}
Then for $\a\in\ell_1(\Z^2)$ and $\bu\in\ell_2(\bOm)$,  we define the convolution $\a\circledast\bu$ as
\begin{align}\label{ConvolutionMeaning}
\left(\a\circledast\bu\right)[\bk]=\left\la\bmP_N(\a)[\bk-\cdot],\bu\right\ra=\sum_{\bn\in\bOm}\bmP_N\left(\a\right)[\bk-\bn]\bu[\bn],~~~~~\bk\in\bOm.
\end{align}
Using the filters defined as \cref{FilterDefinition}, we define the two dimensional fast (discrete) framelet transform, or the analysis operator (see, e.g., \cite{B.Dong2013}) with $L$ levels of decomposition as
\begin{align}\label{Decomposition}
\bmW\bu=\big\{\bmW_{l,\aal}\bu:(l,\aal)\in\big(\left\{0,\cdots,L-1\right\}\times\BB\big)\cup\left\{(L-1,\0)\right\}\big\}.
\end{align}
where $\BB=\left\{0,\cdots,r\right\}^2\setminus\left\{\0\right\}$ is the framelet band. Then $\bmW$ is a linear operator with the frame coefficients $\bmW_{l,\aal}\bu\in\ell_2(\bOm)$ of $\bu$ at level $l$ and band $\aal$ being defined as
\begin{align*}
\bmW_{l,\aal}\bu=\a_{l,\aal}[-\cdot]\circledast\bu,
\end{align*}
where the discrete convolution $\circledast$ is defined as \cref{ConvolutionMeaning}. The synthesis framelet transform is denoted as $\bmW^T$, the adjoint of $\bmW$. Since we use a tight wavelet frame, we have the following perfect reconstruction formula
\begin{align}\label{PerfectRecon}
\bu=\bmW^T\bmW\bu
\end{align}
for all $\bu\in\ell_2(\bOm)$.

Finally, using this analysis operator $\bmW$, we define $\left\|\lam\cdot\bmW\bu\right\|_1$ by
\begin{align}\label{Weightedl1norm}
\left\|\lam\cdot\bmW\bu\right\|_1=\sum_{\bk\in\bOm}\sum_{l=0}^{L-1}\sum_{\aal\in\BB}\lambda_{l,\aal}[\bk]\left|\big(\bmW_{l,\aal}\bu\big)[\bk]\right|.
\end{align}
In particular, if $\lambda_{l,\aal}[\bk]=1$ for all $(l,\aal,\bk)$, we simply write it as $\left\|\bmW\bu\right\|_1$:
\begin{align}\label{Waveletl1normSimple}
\left\|\bmW\bu\right\|_1=\sum_{\bk\in\bOm}\sum_{l=0}^{L-1}\sum_{\aal\in\BB}\left|\big(\bmW_{l,\aal}\bu\big)[\bk]\right|.
\end{align}
As we shall see later, the parameter $\lam$ in \cref{Weightedl1norm} is related to the decay of wavelet frame coefficients, thereby to control the regularity of $\bu$.

\subsection{Error analysis of discrete model \cref{WaveletFrameModel}}\label{WaveletFrameAnalysis}

Throughout the rest of this paper, we denote $\Ps=\{\psi_{\aal}:\aal\in\BB\}$ as the set of framelets, $\phi$ (or $\psi_{\0}$) as the corresponding refinable function, and $\{\a_{\aal}:\aal\in\BB\cup\left\{\0\right\}\}$ as the associated filters. In this paper, we focus on the tensor product B-spline wavelet frame systems constructed by \cite{A.Ron1997}, and $\phi$ is a tensor product B-spline function. We shall refer to $\psi_{\aal}$ for $\aal\in\BB$ as the tensor product framelets. Note, however, that the analysis can be generalized to the generic MRA based wavelet frame systems whenever the corresponding refinable function is compactly supported with $\int_{\R^2}\phi(\x)\rd\x=1$, Lipschitz continuous, and forms a Riesz basis on the closed shift invariant space $\bigvee\left\{\phi(\cdot-\bk):\bk\in\Z^2\right\}$ \cite{J.F.Cai2011}.

Let $\bu^{\La}$ be a solution to \cref{WaveletFrameModel} with the weighted $\ell_1$-norm of wavelet frame coefficients defined as  \cref{Weightedl1norm}. Recall that $\bOm=\left\{0,1,\cdots,N-1\right\}^2$ and $\La$ is a data set uniformly randomly drawn from all $m$-subsets of $\bOm$. Denote by $\rho=m/|\bOm|$ the density of the pixels available. To draw an error between $\bu^{\La}$ and the underlying image $\bsf$, we assume that the linear operator $\bmA$, and the parameter $\lam:=\left\{\lambda_{l,\aal}[\bk]:l=0,\cdots,L-1,~\aal\in\BB,~\bk\in\bOm\right\}$ satisfy the followings:
\begin{enumerate}
\item[A1.] $\s_{\min}(\bmA)>0$, and both $\s_{\min}(\bmA)$ and $\left\|\bmA\right\|_{\infty}$ are independent of $|\bOm|$, i.e. independent of $N$.
\item[A2.] $\lam=\left\{\lambda_{l,\aal}[\bk]:l=0,\cdots,L-1,~\aal\in\BB,~\bk\in\bOm\right\}$ satisfies $\lambda_{l,\aal}[\bk]=2^{\beta\Upsilon(l,\aal,\bk)}$ where $\beta\geq-1$ and $\Upsilon:\Gamma\to\N$ is a positive integer-valued funciton such that
\begin{align}\label{lambdaCondition}
\max\left\{\Upsilon(l,\aal,\bk):(l,\aal,\bk)\in\Gamma\right\}\leq\f{1}{2}\log_2|\bOm|
\end{align}
with $\Gamma=\left\{(l,\aal,\bk):l=0,\cdots,L-1,~\aal\in\BB,~\bk\in\bOm\right\}$.
\end{enumerate}

\begin{example} Before we continue, we present some examples related to A1.
\begin{enumerate}
\item Image inpainting: in this case, we have $\bmA=\bmI$, i.e. $\s_{\min}\left(\bmI\right)=\left\|\bmI\right\|_{\infty}=1$.
\item Simultaneous Gaussian deblurring and image inpainting: let $\bmA\bu=\a\circledast\bu$ where $\a$ is defined as
\begin{align*}
\a[\bk]=\f{C}{2\pi\s^2}\exp\left\{-\f{\left|\bk\right|^2}{2\s^2}\right\},~~~~~\bk\in\Z^2,
\end{align*}
for $\s>0$ and a normalizing constant $C>0$ is such that $\left\|\a\right\|_{\ell_1(\Z^2)}=1$. Obviously, we have $\left\|\bmA\right\|_{\infty}=\left\|\a\right\|_{\ell_1(\Z^2)}=1$. In addition, by Poisson summation formula (e.g. \cite{Mallat2008}), we have
\begin{align*}
\wh{\a}(\oom)=C\sum_{\oom'\in\Z^2}\exp\left\{-\f{\s^2\left|\oom-2\pi\oom'\right|^2}{2}\right\}\geq C\exp\left\{-\s^2\pi^2\right\},~~~\oom\in[-\pi,\pi]^2.
\end{align*}
Hence, we can set $\s_{\min}\left(\bmA\right)=C\exp\left\{-\s^2\pi^2\right\}>0$, which is independent of $|\bOm|$.
\item Wavelet inpainting: let $\bmA$ be a two dimensional orthonormal wavelet basis transform ($\bmA^T\bmA=\bmA\bmA^T=\bmI$) with filter banks $\left\{\bsh_{\bbe}:\bbe\in\left\{0,1\right\}^2\right\}$ generated by the tensor product of a real-valued univariate conjugate mirror filter $\bsh$:
    \begin{align*}
    \bsh_{\bbe}[\bk]=\bsh_{\beta_1}[k_1]\bsh_{\beta_2}[k_2],~~~~~\bk\in\Z^2,~~~\bbe=(\beta_1,\beta_2)\in\left\{0,1\right\}^2
    \end{align*}
    with $\bsh_0=\bsh$ and $\bsh_1[k]=(-1)^{1-k}\bsh[1-k]$. Then $\left\|\bmA\right\|_{\infty}=\left\|\bsh\right\|_{\ell_1(\Z)}^2$ and $\s_{\min}\left(\bmA\right)=\s_{\max}\left(\bmA\right)=1$.
\end{enumerate}
\end{example}

\begin{rmk} Together with the assumption on $\bsf$ that $\left\|\lam\cdot\bmW\bsf\right\|_1<\infty$, A2 is related to the decay of the canonical wavelet frame coefficients $\bmW\bsf$. Notice that the decay of $\bmW\bsf$ is further related to the regularity of $\bsf$, which will be clearer in the context of asymptotic function approximation. See \cref{FunctionApproximation} for details.
\end{rmk}

Under the above assumptions, we note that the underlying true image $\bsf$ satisfies the constraints in \cref{WaveletFrameModel} with $\left\|\lam\cdot\bmW\bsf\right\|_1<\infty$. In addition, since $\bu^{\La}$ is a solution to \cref{WaveletFrameModel}, its weighted $\ell_1$-norm $\left\|\lam\cdot\bmW\bu^{\La}\right\|_1$ attains the minimum subject to the constraints, and in particular, $\left\|\lam\cdot\bmW\bu^{\La}\right\|_1\leq\left\|\lam\cdot\bmW\bsf\right\|_1$. This leads us to consider the following set
\begin{align}\label{HypothesisSpace}
\msM=\left\{\bu\in\ell_{\infty}(\bOm):\left\|\lam\cdot\bmW\bu\right\|_1\leq\left\|\lam\cdot\bmW\bsf\right\|_1,~\f{1}{|\La|}\sum_{\bk\in\La}\left|(\bmA\bu)[\bk]-\bsg[\bk]\right|^2\leq\eta^2,~\bu\in[0,M]^{\bOm}\right\}
\end{align}
as an involved hypothesis space. In \cref{HypothesisSpace}, $M>0$ is a positive constant related to the boundedness of each pixel value, $\eta>0$ is a fixed positive constant related to the bound of measurement error, and the weighted $\ell_1$-norm of wavelet frame coefficients $\left\|\lam\cdot\bmW\bu\right\|_1$ is defined as \cref{Weightedl1norm}.

Our error estimate is based on the characterization of the capacity of the hypothesis space \cite{J.F.Cai2011}. Among such capacities including VC dimension \cite{Vapnik1998}, $V_{\gamma}$-dimension and $P_{\gamma}$-dimension \cite{N.Alon1997}, Rademacher complexities \cite{P.L.Bartlett2006,Koltchinskii2001}, and covering number \cite{F.Cucker2002}, we choose the covering number because it is the most convenient for metric spaces \cite{J.F.Cai2011}.

\begin{definition} Let $\msM\subseteq\R^{\bOm}$ and $r>0$ be given. The covering number $\mN(\msM,r)$ is defined as
\begin{align*}
\mN\left(\msM,r\right)=\inf\left\{K\in\N:\exists~\bu_1,\cdots,\bu_K\in\msM~~\text{s.t.}~~\msM\subseteq\bigcup_{j=1}^{K}\left\{\bu\in\msM:\left\|\bu-\bu_j\right\|_{\ell_{\infty}(\bOm)}\leq r\right\}\right\},
\end{align*}
i.e. the minimal number of $\ell_{\infty}$ balls with radius $r$ in $\msM$ that covers $\msM$.
\end{definition}

With the above notation, we can present the first relation between the solution $\bu^{\La}$ to \cref{WaveletFrameModel} and the underlying image $\bsf$. Following the idea of statistical learning theory \cite{Vapnik1998}, we first establish the probability of the following event
\begin{align*}
\f{1}{|\bOm|}\left\|\bu^{\La}-\bsf\right\|_{\ell_2(\bOm)}^2\leq\ep+\f{16}{3}\s_{\min}^{-2}(\bmA)\eta^2
\end{align*}
for an arbitrary $\ep>0$ in terms of the covering number with the radius related to $\ep$. The proof is postponed to \cref{ProofTh1}

\begin{theorem}\label{Th1} Let $\msM$ be defined as in \cref{HypothesisSpace} and $\bu^{\La}$ be a solution to \cref{WaveletFrameModel}. Then for an arbitrary $\ep>0$, the following inequality
\begin{align*}
\PP\left\{\f{1}{|\bOm|}\left\|\bu^{\La}-\bsf\right\|_{\ell_2(\bOm)}^2\leq\ep+\f{16}{3}\s_{\min}^{-2}(\bmA)\eta^2\right\}\geq1-\mN\left(\msM,\f{\s_{\min}^2(\bmA)\ep}{12\left\|\bmA\right\|_{\infty}^2M}\right)\exp\left\{-\f{3m\s_{\min}^2(\bmA)\ep}{256\left\|\bmA\right\|_{\infty}^2M^2}\right\}
\end{align*}
holds for an arbitrary $m$, where $m=|\La|$ denotes the number of samples.
\end{theorem}

\begin{rmk} Notice that the proof of \cref{Th1} is completed if we estimate the probability of the event
\begin{align*}
\f{1}{|\bOm|}\left\|\bmA\bu^{\La}-\bmA\bsf\right\|_{\ell_2(\bOm)}^2\leq\s_{\min}^2\left(\bmA\right)\ep+\f{16}{3}\eta^2.
\end{align*}
In this case, we need to use the uniform boundedness constraint ($\bu[\bk]\in[0,M]$ for all $\bk\in\bOm$), from which we need an assumption on $\left\|\bmA\right\|_{\infty}$. See \cref{ProofTh1} for details.
\end{rmk}

By \cref{Th1}, our error estimate will be completed if we estimate the covering number $\mN\left(\msM,r\right)$. Intuitively, since $\msM$ lies in an $\ell_{\infty}$ ball $\big\{\bu\in\ell_{\infty}(\bOm):\left\|\bu\right\|_{\ell_{\infty}(\bOm)}\leq M\big\}$, it is obvious that
\begin{align}\label{CoveringNumberRoughEstimate}
\mN(\msM,r)\leq\left(\f{2M}{r}\right)^{|\bOm|}
\end{align}
as given in \cite{F.Cucker2002}. However, since the above estimation \cref{CoveringNumberRoughEstimate} is not tight enough to derive an error estimate, we need to find a tighter upper bound for $\mN(\msM,r)$ to obtain a reasonable error bound. To do this, we adopt the idea given in \cite{J.F.Cai2011}; using the discrete total variation, we connect the wavelet frame coefficients $\bmW\bu$ with the regularity of $\bu$. Specifically, we define the discrete gradient of $\bu\in\R^{\bOm}$ as
\begin{align}\label{DiscreteGradient}
\Na\bu=\left\{\bu[\bk+\e_j]-\bu[\bk]:j=1,2,~\bk\in\bOm,~\&~\bk+\e_j\in\bOm\right\},
\end{align}
where $\e_1=(1,0)$ and $\e_2=(0,1)$. The discrete total variation $\left\|\Na\bu\right\|_1$ is defined as
\begin{align}\label{DiscreteTV}
\left\|\Na\bu\right\|_1:=\sum_{j=1}^2\sum_{\stackrel{\bk\in\bOm}{\bk+\e_j\in\bOm}}\left|\bu[\bk+\e_j]-\bu[\bk]\right|.
\end{align}
Under this setting, we present \cref{Lemma1} whose proof is postponed to \cref{ProofLemma1}.

\begin{lemma}\label{Lemma1}
Let $\bmW$ be a tensor product B-spline wavelet frame transform defined as \cref{Decomposition}. Then there exists a constant $C_{\bmW}\geq1$ which is independent of $|\bOm|$ and $L$ such that
\begin{align}\label{Lemma1Result}
\left\|\na\bu\right\|_1\leq C_{\bmW}\left\|\bmW\bu\right\|_1
\end{align}
for all $\bu\in\R^{\bOm}$. In \cref{Lemma1Result}, $\left\|\Na\bu\right\|_1$ is defined as \cref{DiscreteTV}, and $\left\|\bmW\bu\right\|_1$ is defined as \cref{Waveletl1normSimple}.
\end{lemma}

\cref{Lemma1} tells us that, given that $\lam=\left\{\lambda_{l,\aal}[\bk]\right\}$ satisfies A2, for $\bu\in\msM$, we have
\begin{align*}
\left\|\Na\bu\right\|_1&\leq C_{\bmW}\left\|\bmW\bu\right\|_1=C_{\bmW}\sum_{\bk\in\bOm}\sum_{l=0}^{L-1}\sum_{\aal\in\BB}\left|\left(\bmW_{l,\aal}\bu\right)[\bk]\right|\\
&\leq C_{\bmW}2^{\max\left\{-\beta,0\right\}\max\left\{\Upsilon(l,\aal,\bk):(l,\aal,\bk)\in\Gamma\right\}}\sum_{\bk\in\bOm}\sum_{l=0}^{L-1}\sum_{\aal\in\BB}2^{\beta\Upsilon(l,\aal,\bk)}\left|\left(\bmW_{l,\aal}\bu\right)[\bk]\right|\\
&\leq C_{\bmW}|\bOm|^{\f{\max\{1-\beta,0\}}{2}}\left\|\lam\cdot\bmW\bsf\right\|_1.
\end{align*}
Hence, we can further relax the set $\msM$ into
\begin{align}\label{RelaxedFeasibleSet}
\wt{\msM}=\left\{\bu\in\ell_{\infty}(\bOm):\left\|\na\bu\right\|_1\leq C_{\bmW}|\bOm|^{\f{\max\{1-\beta,0\}}{2}}\left\|\lam\cdot\bmW\bsf\right\|_1,~\left\|\bu\right\|_{\ell_{\infty}(\bOm)}\leq M\right\}.
\end{align}
Then with the set $\wt{\msM}$ defined above, we present \cref{Th4} for the estimation of covering number. Similar to \cite[Theorem 2.4]{J.F.Cai2011}, our estimate is also based on the quantized total variation minimization (e.g. \cite{A.Chambolle[2021]copyright2021,A.Chambolle2009}) and the max-flow min-cut theorem \cite{Diestel2018} to introduce a new upper bound of the covering number.

\begin{theorem}\label{Th4} Let $\bmW$ be a tensor product B-spline wavelet frame transform defined as \cref{Decomposition}. Assume that A2 is satisfied with $-1\leq\beta\leq 1$ and $\left\|\lam\cdot\bmW\bsf\right\|_1\leq C_{\bsf}|\bOm|^{1/2}$. Let $\msM$ be defined as \cref{HypothesisSpace}. For $r\geq|\bOm|^{-a}$ with $a\geq1$, we have
\begin{align*}
\ln\mN(\msM,r)\leq\f{C_a|\bOm|^{\f{\max\{1-\beta,1\}}{2}}}{r}\log_2|\bOm|,
\end{align*}
where $C_a=20M(4a+\max\{1-\beta,1\})C_{\bmW}C_{\bsf}(1+C_{\bmW}C_{\bsf})$.
\end{theorem}

\begin{pf} See \cref{ProofTh4}.
\end{pf}

With the aid of \cref{Th1,Th4}, we can establish the explicit formulation of \cref{Goal} in \cref{Th2}. Briefly speaking, for a fixed $\rho$, as long as the resolution of an image is sufficiently high, the wavelet frame based restoration model \cref{WaveletFrameModel} has a good chance to restore the original image within the measurement error. In addition, for a fixed $|\bOm|$, the error gets smaller as $\rho\nearrow1$, which coincides with our common sense. Finally, \cref{Th2} takes the results in \cite{J.F.Cai2011} as a special case because $\left\|\bmI\right\|_{\infty}=\s_{\min}(\bmI)=1$.

\begin{theorem}\label{Th2} Let $\bmW$ be a tensor product B-spline wavelet frame transform. Assume that $\bmA$ satisfies A1, $\lam$ satisfies A2 with $-1<\beta<1$, and $\bsf\in[0,M]^{\bOm}$ satisfies $\left\|\lam\cdot\bmW\bsf\right\|_1\leq C_{\bsf}|\bOm|^{1/2}$. Let $\bu^{\La}$ be a solution to \cref{WaveletFrameModel}. Then the following inequality
\begin{align}\label{MainResult:Discrete}
\f{1}{|\bOm|}\left\|\bu^{\La}-\bsf\right\|_{\ell_2(\bOm)}^2\leq\wt{c}\rho^{-\f{1}{2}}|\bOm|^{-\f{\min\{1+\beta,1\}}{4}}\left(\log_2|\bOm|\right)^{3/2}+\f{16}{3}\s_{\min}^{-2}(\bmA)\eta^2
\end{align}
holds with probability at least $1-|\bOm|^{-1}$. In \cref{MainResult:Discrete}, $\wt{c}$ is defined as
\begin{align*}
\wt{c}=\f{64\left\|\bmA\right\|_{\infty}^2M^2}{3\s_{\min}^2(\bmA)}\left[4+3\sqrt{5(4a+\max\{1-\beta,1\})C_{\bmW}C_{\bsf}(1+C_{\bmW}C_{\bsf})}\right],
\end{align*}
for some fixed $a\geq1$. In words, $\wt{c}$ is independent of $|\bOm|$, $L$, $\rho$, and $\eta$.
\end{theorem}

\begin{pf} First of all, for any $\ep>0$, we can choose $a\geq1$ such that $\ep\geq\f{12\left\|\bmA\right\|_{\infty}^2M}{\s_{\min}^2(\bmA)}|\bOm|^{-a}$. With $r=\f{\s_{\min}^2(\bmA)\ep}{12\left\|\bmA\right\|_{\infty}^2M}\geq|\bOm|^{-a}$, \cref{Th1,Th4} give us that the inequality
\begin{align*}
\f{1}{|\bOm|}\left\|\bu^{\La}-\bsf\right\|_{\ell_2(\bOm)}^2\leq\ep+\f{16}{3}\s_{\min}^{-2}(\bmA)\eta^2
\end{align*}
holds with probability at least
\begin{align*}
1&-\mN\left(\msM,\f{\s_{\min}^2(\bmA)\ep}{12\left\|\bmA\right\|_{\infty}^2M}\right)\exp\left\{-\f{3m\s_{\min}^2(\bmA)\ep}{256\left\|\bmA\right\|_{\infty}^2M^2}\right\}\\
&\geq1-\exp\left\{\f{240\left\|\bmA\right\|_{\infty}^2M^2(4a+b)C_{\bmW}C_{\bsf}(1+C_{\bmW}C_{\bsf})|\Om|^{b/2}\log_2|\bOm|}{\s_{\min}^2(\bmA)\ep}-\f{3m\s_{\min}^2(\bmA)\ep}{256\left\|\bmA\right\|_{\infty}^2M^2}\right\}.
\end{align*}
In the above, we use $b=\max\{1-\beta,1\}$ for notational simplicity. Hence, if we choose a special $\ep^*$ to be the unique positive solution to the following equation
\begin{align}\label{Eq:epsilon}
\f{240\left\|\bmA\right\|_{\infty}^2M^2(4a+b)C_{\bmW}C_{\bsf}(1+C_{\bmW}C_{\bsf})|\Om|^{b/2}\log_2|\bOm|}{\s_{\min}^2(\bmA)\ep}-\f{3m\s_{\min}^2(\bmA)\ep}{256\left\|\bmA\right\|_{\infty}^2M^2}=\ln\f{1}{|\bOm|},
\end{align}
we have
\begin{align*}
\f{1}{|\bOm|}\left\|\bu^{\La}-\bsf\right\|_{\ell_2(\bOm)}^2\leq\ep^*+\f{16}{3}\s_{\min}^{-2}(\bmA)\eta^2
\end{align*}
with probability at least $1-|\bOm|^{-1}$. Indeed, solving \cref{Eq:epsilon} yields \cref{MainResult:Discrete}:
\begin{align*}
\ep^*&=\f{128\left\|\bmA\right\|_{\infty}^2M^2}{3m\s_{\min}^2(\bmA)}\left[\ln|\bOm|+\sqrt{\ln^2|\bOm|+\f{45}{4}m(4a+b)C_{\bmW}C_{\bsf}(1+C_{\bmW}C_{\bsf})|\bOm|^{b/2}\log_2|\bOm|}\right]\\
&\leq\f{64\left\|\bmA\right\|_{\infty}^2M^2}{3m\s_{\min}^2(\bmA)}\left[4\ln|\bOm|+3\sqrt{5m(4a+b)C_{\bmW}C_{\bsf}(1+C_{\bmW}C_{\bsf})(1+C_{\bmW}C_{\bsf})|\bOm|^{b/2}\log_2|\bOm|}\right]\\
&\leq\f{64\left\|\bmA\right\|_{\infty}^2M^2}{3\s_{\min}^2(\bmA)}\left[4+3\sqrt{5(4a+b)C_{\bmW}C_{\bsf}(1+C_{\bmW}C_{\bsf})}\right]\rho^{-\f{1}{2}}|\bOm|^{\f{b-2}{4}}\left(\log_2|\bOm|\right)^{3/2},
\end{align*}
as $b-2=-\min\{1+\beta,1\}$. In addition, from the definition of $\ep^*$, we have
\begin{align*}
\ep^*&\geq\f{128\left\|\bmA\right\|_{\infty}^2M^2}{\s_{\min}^2(\bmA)\sqrt{m}}\sqrt{C_{\bmW}C_{\bsf}(1+C_{\bmW}C_{\bsf})|\bOm|^{b/2}\log_2|\bOm|}\\
&\geq\f{16\left\|\bmA\right\|_{\infty}^2M}{\s_{\min}^2(\bmA)}\sqrt{C_{\bmW}C_{\bsf}(1+C_{\bmW}C_{\bsf})}\rho^{-\f{1}{2}}\sqrt{\log_2|\bOm|}|\bOm|^{\f{b-2}{4}}\geq\f{12\left\|\bmA\right\|_{\infty}^2M}{\s_{\min}^2(\bmA)}|\bOm|^{-a}.
\end{align*}
In the final inequality, we use the fact that for $\beta\in(-1,1)$, we have $b-2=-\min\{1+\beta,1\}\in(-1,0)$ while $a\geq1$. This completes the proof.\qquad$\square$
\end{pf}

\begin{rmk} From the proof of \cref{Th2}, we can see that the constant $a\geq1$ is to apply \cref{Th4} to \cref{Th1}. Notice that the lower bound of the probability of
\begin{align*}
\f{1}{|\bOm|}\left\|\bu^{\La}-\bsf\right\|_{\ell_2(\bOm)}^2\leq\ep+\f{16}{3}\s_{\min}^{-2}\left(\bmA\right)\eta^2
\end{align*}
decreases as $\ep>0$ decreases. Meanwhile, we estimate the critical $\ep^*>0$ (i.e. the positive solution to \cref{Eq:epsilon}) by fixing the probability to be no smaller than $1-|\bOm|^{-1}$, which in turn increases the upper bound (or the confidence interval) of $\ep^*>0$ as reflected by the constant $a\geq1$. Nevertheless, we forgo the effect of $a\geq1$ in \cref{MainResult:Discrete} to the error estimate in order not to dilute the focus of this paper. In fact, $a\geq1$ becomes fixed once we fix the admissible lower bound of $\ep>0$ in \cref{Th1} to apply \cref{Th4}.
\end{rmk}

To conclude, we present some numerical results to demonstrate the theoretical error bound in \cref{MainResult:Discrete} and the empirical restoration error under different settings of $|\bOm|$ and $\rho$. Specifically, we fix $N=2^J$ and $\eta=0$, and we consider the following noise-free case
\begin{align}\label{Model:Goal_NoiseFree}
\min_{\bu\in[0,M]^{\bOm}}~\left\|\lam\cdot\bmW\bu\right\|_1~~~\text{subject to}~~~\left(\bmA\bu\right)[\bk]=\bsg[\bk],~~\bk\in\La,
\end{align}
so that \cref{MainResult:Discrete} takes the form of
\begin{align}\label{MainResult:NoiseFree}
\f{1}{2^{2J}}\left\|\bu^{\La}-\bsf\right\|_{\ell_2(\bOm)}^2\leq\wt{c}\rho^{-1/2}J^{3/2}2^{-J\min\left\{\f{1+\beta}{2},\f{1}{2}\right\}}
\end{align}
with probability at least $1-2^{-2J}$.

\begin{figure}[tp]
\centering
\subfloat[]{\label{FrameInpaintingResultsSampleDensity}\includegraphics[width=8.1cm]{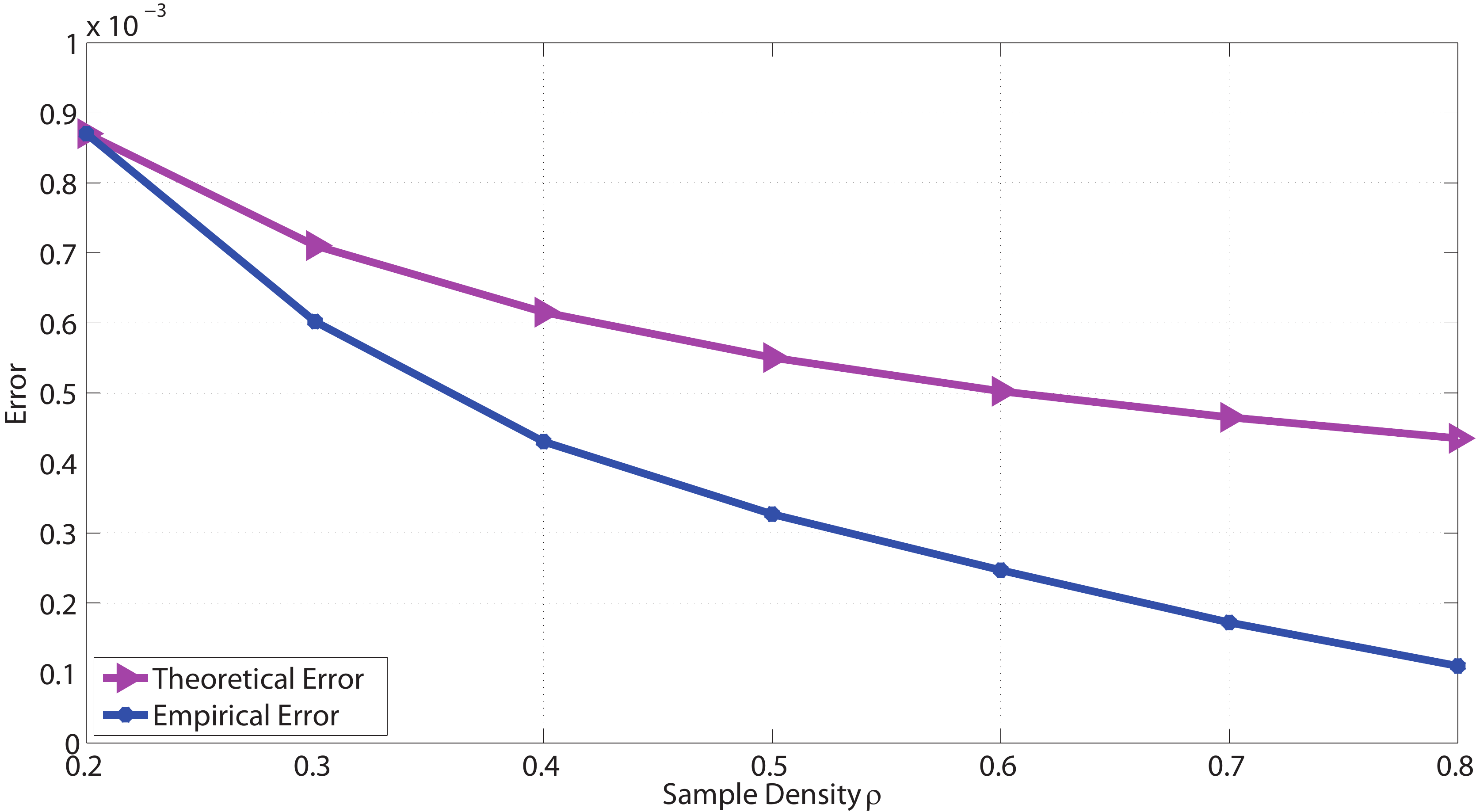}}\hspace{0.001cm}
\subfloat[]{\label{FrameInpaintingResultsResolution}\includegraphics[width=8.1cm]{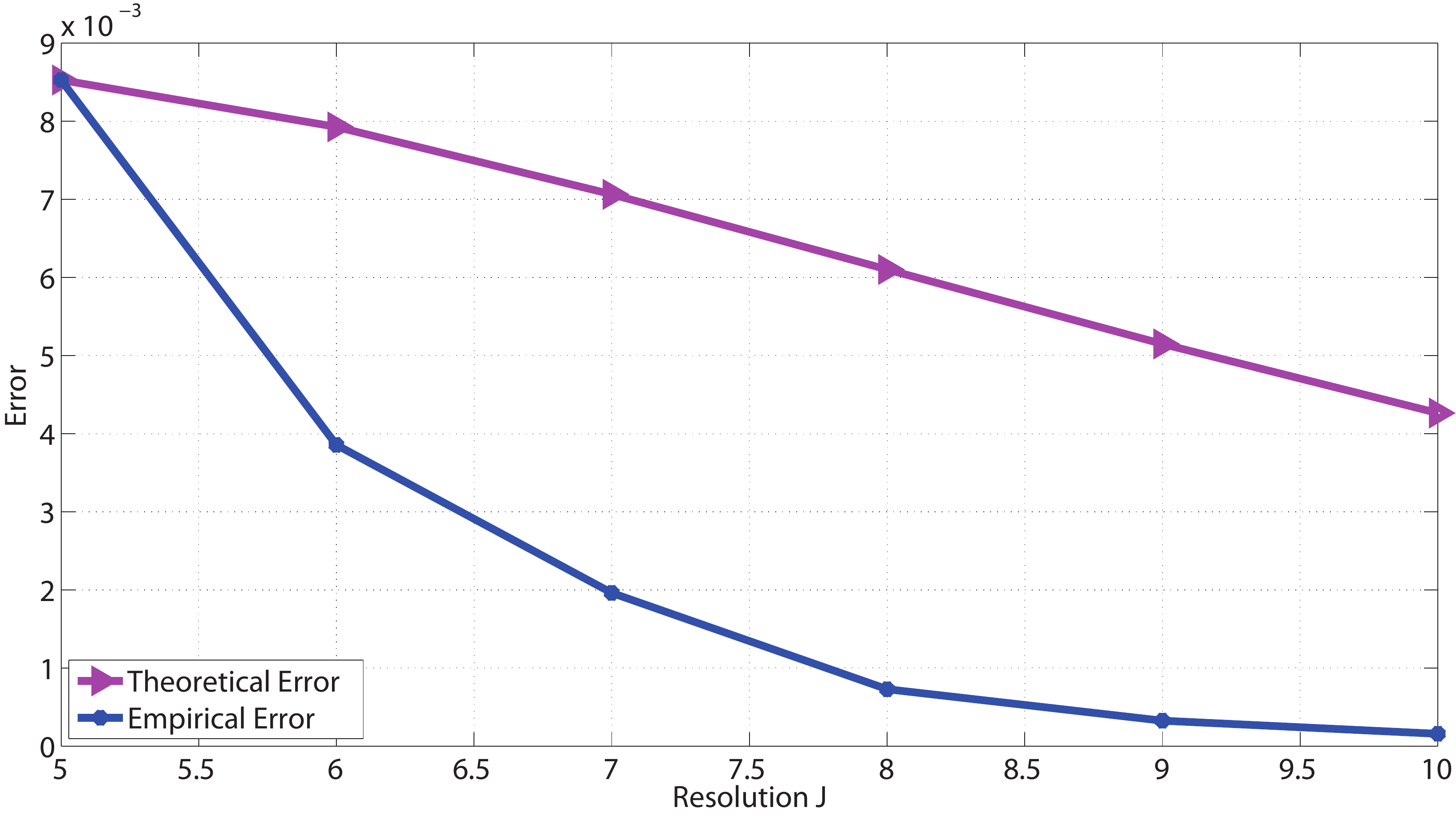}}\hspace{0.001cm}
\caption{Simulation results for the image inpainting. \cref{FrameInpaintingResultsSampleDensity} describes the simulation results when $N=512$ is fixed, and \cref{FrameInpaintingResultsResolution} depicts the simulation results when $\rho=0.5$ is fixed.}\label{FrameInpaintingResultsAnalysis}
\end{figure}

\begin{figure}[tp]
\centering
\subfloat[]{\label{FrameDebInpResultsSampleDensity}\includegraphics[width=8.1cm]{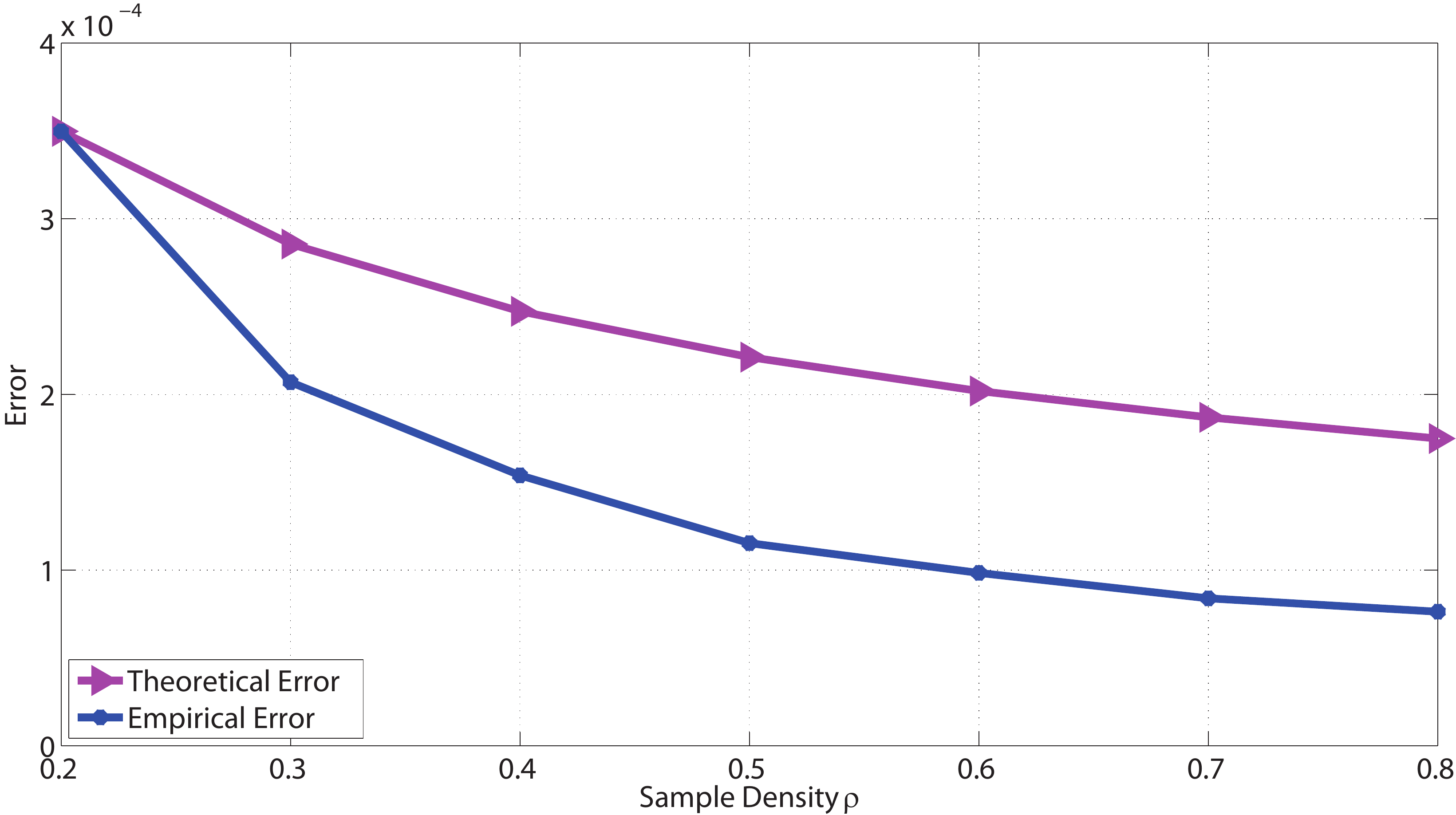}}\hspace{0.001cm}
\subfloat[]{\label{FrameDebInpResultsResolution}\includegraphics[width=8.1cm]{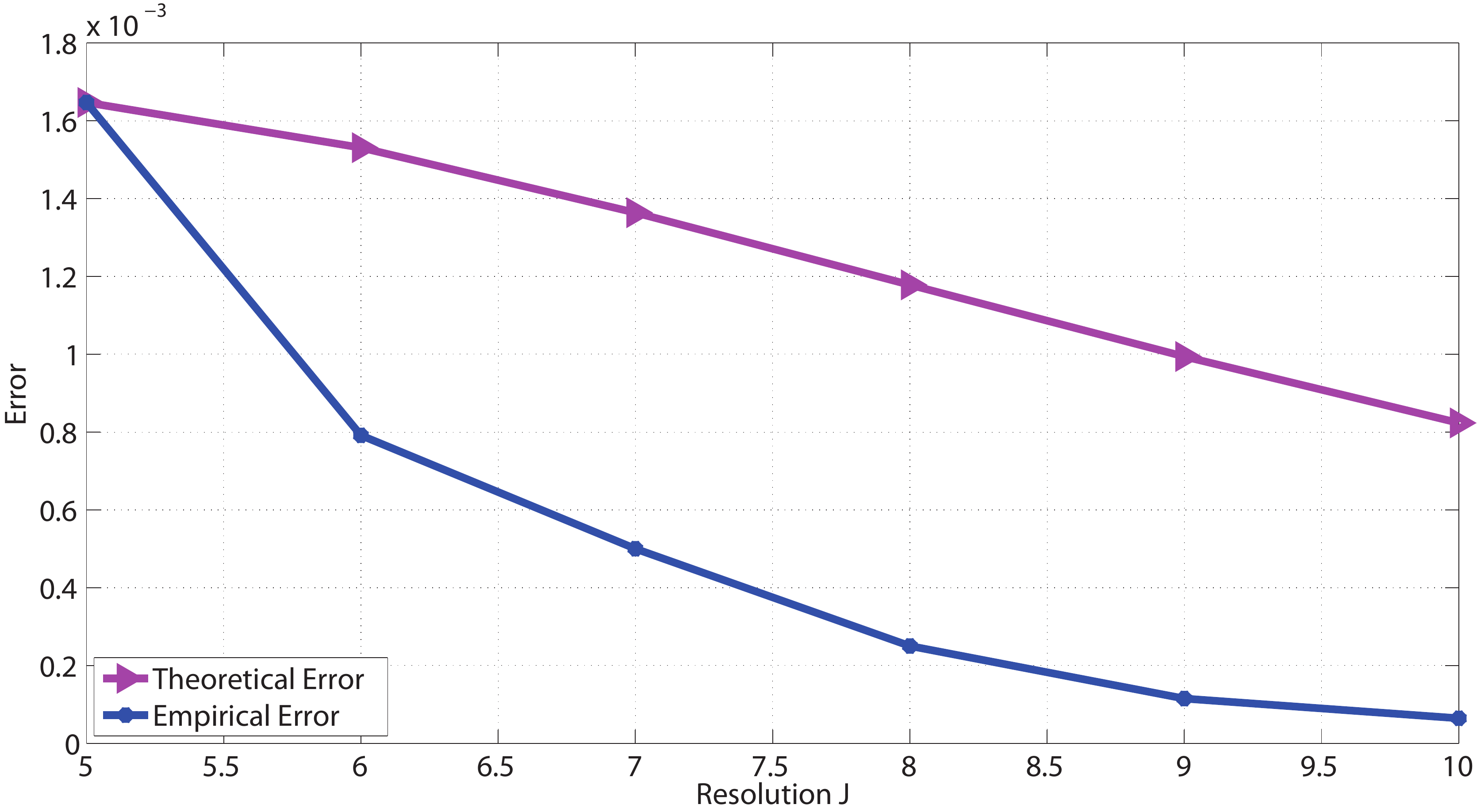}}\hspace{0.001cm}
\caption{Simulation results for the simultaneous deblurring and inpainting. \cref{FrameDebInpResultsSampleDensity} describes the simulation results when $N=512$ is fixed, and \cref{FrameDebInpResultsResolution} depicts the simulation results when $\rho=0.5$ is fixed.}\label{FrameDebInpResultsAnalysis}
\end{figure}

In this simulation, two types of $\bmA$ are used; one is the identity operator $\bmI$ (i.e. the image inpainting) and the other is a convolution operator with a low-pass filter generated by the Matlab built-in function ``fspecial(`gaussian',$9$,$1$)'' (i.e. the simultaneous deblurring and inpainting). More precisely, we generate a piecewise linear phantom as an underlying true image $\bsf$, and we implement the simulation for each $\bmA$ as follows. To see the behavior of error with respect to the sample density $\rho$, we fix $N=512$ (i.e. $J=9$), and for each $\rho\in\left\{0.2,0.3,0.4,0.5,0.6,0.7,0.8\right\}$, we test \cref{Model:Goal_NoiseFree} with $100$ realizations of $\La$. To see the behavior of error with respect to $|\bOm|$, we fix $\rho=0.5$, and for each $J\in\left\{5,6,7,8,9,10\right\}$ (or the resolution $J$), we again test \cref{Model:Goal_NoiseFree} with $100$ realizations of $\La$. In any case, we choose $\bmW$ to be the tensor product piecewise linear B-spline wavelet frame with $1$ level of decomposition, and \cref{Model:Goal_NoiseFree} is solved by the split Bregman algorithm \cite{J.F.Cai2009/10}. Finally, we choose the largest empirical error to compare with the smallest theoretical error in \cref{MainResult:NoiseFree} by setting $\beta=0$ in \cref{MainResult:NoiseFree}. In particular, it is in general difficult to determine the constants explicitly, we calculate the above error by assuming that the equality holds in the worst case, i.e. in the lowest sample density, or the lowest resolution case. The results are presented in \cref{FrameInpaintingResultsAnalysis,FrameDebInpResultsAnalysis} respectively. Specifically, \cref{FrameInpaintingResultsSampleDensity,FrameDebInpResultsSampleDensity} show the results when $N=512$ is fixed and $\rho$ is varying, and \cref{FrameInpaintingResultsResolution,FrameDebInpResultsResolution} demonstrate when $\rho=0.5$ is fixed and $N$ is varying. We can easily see that, in any case, the empirical restoration error does not exceed the theoretical error in \cref{MainResult:NoiseFree}, which empirically demonstrates that \cref{Th2} provides a reasonable upper bound for the restoration error with high probability.

\subsection{Asymptotic function approximation via MRA}\label{FunctionApproximation}

In this section, we apply \cref{Th2} to further connect with the function approximation given that the underlying true image is obtained by the samples of a function. All functions we consider are defined on $\Om=[0,1)^2$, and we consider the $2^J\times 2^J$ Cartesian grid $\bOm$ defined as
\begin{align*}
\bOm=\left\{0,1,\cdots,2^J-1\right\}^2\simeq2^{-J}\Z^2\cap\Om.
\end{align*}
In other words, we implicitly identify the $2^J\times2^J$ grid $\bOm$ with the $2^J\times2^J$ mesh of $\Om$. For the consistency with the discrete image cases, we identify $L_p(\Om)$ with the space of $L_p$ functions on $\R^2$ with fundamental period on each variable to be $1$. For the wavelet frame decomposition (e.g. \cite{I.Daubechies2003,B.Dong2013}) of $f\in L_2(\Om)$, we write
\begin{align*}
f=\sum_{l\in\Z}\sum_{\bk\in\bOm}\sum_{\aal\in\BB}\la f,\psi_{\aal,l,\bk}\ra\psi_{\aal,l,\bk}
\end{align*}
by identifying $\psi_{\aal,l,\bk}$ with its periodized version
\begin{align*}
\psi_{\aal,l,\bk}^{\mathrm{per}}=\sum_{\bk'\in\Z^2}\psi_{\aal,l,\bk}(\cdot+2^{\max\{l,J\}}\bk')
\end{align*}
with a slight abuse of notation.

In the literature of wavelet frame, the discrete image $\bsf$ is interpreted as discrete sampling of an underlying function $f\in L_2(\Om)$ via the inner product with the corresponding refinable function $\phi$. More precisely, let $\phi_{J,\bk}=2^J\phi(2^J\cdot-\bk)$. Then we can write $\bsf$ as
\begin{align}\label{Samples}
\bsf[\bk]=2^J\left\la f,\phi_{J,\bk}\right\ra.
\end{align}
With this $\bsf$, we implicitly use the following interpolated function
\begin{align}\label{InterpFunction}
f_J=\sum_{\bk\in\bOm}\bsf[\bk]\phi(2^J\cdot-\bk)=\sum_{\bk\in\bOm}\left\la f,\phi_{J,\bk}\right\ra\phi_{J,\bk}
\end{align}
to approximate $f$.

There are extensive studies on the approximation order of $f_{J}$ in \cref{InterpFunction} to the original $f$, and most of them are based on certain properties of the refinable function $\phi$. Relevant works include the wavelet frame based scattered data restoration \cite{J.Yang2017} and the wavelet frame based image denoising/inpainting \cite{B.Dong2021}. Briefly speaking, these works commonly assume the Strang-Fix condition on $\phi$ (e.g. \cite{I.Daubechies2003,M.J.Johnson2009}) to approximate the underlying function in the Sobolev space of a certain order. Meanwhile, the underlying function that we concern in this paper does not satisfy certain order of regularity. Instead, we only impose some decay condition of the wavelet frame coefficients for the analysis of the approximation order. More precisely, we assume that there is $\beta\geq-1$ such that
\begin{align}\label{DecayCondition}
C_f:=\sum_{l=J-L}^{\infty}2^{\beta l}\sum_{\bk\in\bOm}\sum_{\aal\in\BB}\big|\la f,\psi_{\aal,l,\bk}\ra\big|<\infty.
\end{align}
Note that the above decay condition links to the regularity of the underlying function $f$ when the wavelet frame system satisfies some mild conditions \cite{L.Borup2004,J.F.Cai2011,B.Han2009}.

Under the decay condition given in \cref{DecayCondition}, we present the following approximation of underlying function in \cref{Th3}. Briefly speaking, as long as the image resolution is sufficiently large, the approximated function constructed from the restored image gives a good approximation of the underlying function where the underlying true image comes from.

\begin{theorem}\label{Th3} Let $\bmW$ be a tensor product B-spline wavelet frame transform. Assume that A1 is satisfied, and $\lam=\left\{\lambda_{l,\aal}[\bk]:l=0,\cdots,L-1,~\aal\in\BB,~\bk\in\bOm\right\}$ is defined as
\begin{align}\label{lambdaDef}
\lambda_{l,\aal}[\bk]=2^{\beta(J-l-1)},~~~~~l=0,\cdots,L-1,~\aal\in\BB,~\bk\in\bOm,
\end{align}
and the underlying function $f\in L_2(\Om)$ satisfies \cref{DecayCondition} with $-1<\beta<1$. Let $\bu^{\La}$ be a solution to \cref{WaveletFrameModel} with $\bsg$ in \cref{Linear_IP_Incomplete} generated by $\bsf$ in \cref{Samples}. If $f(\x)\in[0,M]$ a.e. $\x\in\Om$, we have the followings:
\begin{enumerate}
\item The inequality
\begin{align}\label{CorollaryResult1}
\f{1}{2^{2J}}\left\|\bu^{\La}-\bsf\right\|_{\ell_2(\bOm)}^2\leq\wt{c}\rho^{-\f{1}{2}}J^{3/2}2^{-J\min\left\{\f{1+\beta}{2},\f{1}{2}\right\}}+\f{16}{3}\s_{\min}^{-2}(\bmA)\eta^2
\end{align}
holds with probability at least $1-2^{-2J}$, where $\wt{c}$ is a constant independent of $J$ (i.e. independent of $|\bOm|$), $L$, $\rho$, and $\eta$.
\item Let $u_J^{\La}$ be defined as
\begin{align*}
u_J^{\La}=\sum_{\bk\in\bOm}\bu^{\La}[\bk]\phi(2^J\cdot-\bk).
\end{align*}
Then
\begin{align}\label{CorollaryResult2}
\left\|u_J^{\La}-f\right\|_{L_2(\Om)}^2\leq C_1\rho^{-\f{1}{2}}J^{3/2}2^{-J\min\left\{\f{1+\beta}{2},\f{1}{2}\right\}}+C_2\eta^2+C_32^{-(\beta+1)J}
\end{align}
holds with probability at least $1-2^{-2J}$, where $C_1$, $C_2$, and $C_3$ are three positive constants independent of $J$, $L$, $\rho$, and $\eta$.
\end{enumerate}
\end{theorem}

\begin{rmk} Before proving \cref{Th3}, we would like to mention that, whenever $J\in\N$ satisfies
\begin{align*}
J\geq-\f{2\log_2\eta}{\beta+1},
\end{align*}
\cref{CorollaryResult2} can be further written as
\begin{align*}
\left\|u_J^{\La}-f\right\|_{L_2(\Om)}^2\leq C_1\rho^{-\f{1}{2}}J^{3/2}2^{-J\min\left\{\f{1+\beta}{2},\f{1}{2}\right\}}+\left(C_2+C_3\right)\eta^2.
\end{align*}
In other words, whenever the mesh is sufficiently dense (i.e. $J$ is sufficiently large), the $L_2$ distance between the interpolation of the restored image and the original unknown function is bounded by the restoration error of the discrete image restoration problem \cref{WaveletFrameModel} only.
\end{rmk}

The rest is devoted to the proof of \cref{Th3}. To begin with, we introduce the following lemma on the Bessel property of $\left\{\phi(2^J\cdot-\bk):\bk\in\bOm\right\}$.

\begin{lemma}\label{Lemma3} Let $\phi$ be a tensor product B-spline of order $n$. For $J\in\N$, we have the followings.
\begin{enumerate}
\item For $u\in L_2(\Om)$, we have
\begin{align}\label{Bessel}
\sum_{\bk\in\bOm}\left|\left\la u,\phi(2^J\cdot-\bk)\right\ra\right|^2\leq\f{n^2}{2^{2J}}\left\|u\right\|_{L_2(\Om)}^2.
\end{align}
\item For $\bu\in\ell_2(\bOm)$, we have
\begin{align}\label{AdjointBessel}
\left\|\sum_{\bk\in\bOm}\bu[\bk]\phi(2^J\cdot-\bk)\right\|_{L_2(\Om)}^2\leq\f{n^2}{2^{2J}}\left\|\bu\right\|_{\ell_2(\bOm)}^2.
\end{align}
\end{enumerate}
\end{lemma}

\begin{pf} Notice that $0\leq\phi\leq1$, and the direct computation shows that $\left\|\phi\right\|_{L_2(\R^2)}^2\leq1$, and $\left\|\phi(2^J\cdot-\bk)\right\|_{L_2(\R^2)}^2\leq2^{-2J}$. Then by the Schwarz inequality, we have
\begin{align*}
\left|\left\la u,\phi(2^J\cdot-\bk)\right\ra\right|\leq2^{-J}\left(\int_{\Om}\left|u(\x)\right|^21_{\bS_{\bk}}(\x)\rd\x\right)^{1/2}
\end{align*}
where $\bS_{\bk}=\su\left(\phi(2^J\cdot-\bk)\right)$. Since $\sum 1_{\bS_{\bk}}=n^2$, we have
\begin{align*}
\sum_{\bk\in\bOm}\left|\left\la u,\phi(2^J\cdot-\bk)\right\ra\right|^2\leq2^{-2J}\int_{\Om}\left|u(\x)\right|^2\left(\sum_{\bk\in\bOm}1_{\bS_{\bk}}(\x)\right)\rd\x=\f{n^2}{2^{2J}}\left\|u\right\|_{L_2(\Om)}^2,
\end{align*}
which proves \cref{Bessel}.

For \cref{AdjointBessel}, let $v\in L_2(\Om)$. We have
\begin{align*}
\left\la\sum_{\bk\in\bOm}\bu[\bk]\phi(2^J\cdot-\bk),v\right\ra=\sum_{\bk\in\bOm}\bu[\bk]\left\la\phi(2^J\cdot-\bk),v\right\ra.
\end{align*}
By the Schwarz inequality, we have
\begin{align*}
\sum_{\bk\in\bOm}\left|\bu[\bk]\right|\left|\left\la\phi(2^J\cdot-\bk),v\right\ra\right|\leq\left(\sum_{\bk\in\bOm}\left|\bu[\bk]\right|^2\right)^{1/2}\left(\sum_{\bk\in\bOm}\left|\left\la\phi(2^J\cdot-\bk),v\right\ra\right|^2\right)^{1/2}.
\end{align*}
By \cref{Bessel}, we further have
\begin{align*}
\left|\left\la\sum_{\bk\in\bOm}\bu[\bk]\phi(2^J\cdot-\bk),v\right\ra\right|\leq\left(\sum_{\bk\in\bOm}\left|\bu[\bk]\right|^2\right)^{1/2}\f{n}{2^J}\left\|v\right\|_{L_2(\Om)}.
\end{align*}
Since $v\in L_2(\Om)$ is arbitrary, we have \cref{AdjointBessel} by the converse of H\"older's inequality with $p=q=2$ (e.g. \cite{Folland1999}). This completes the proof.\qquad$\square$
\end{pf}

\begin{poth3} Notice that in this case, we have $|\bOm|=2^{2J}$. In addition, for $\lam$ defined as \cref{lambdaDef}, we have
\begin{align*}
\Upsilon(l,\aal,\bk)=J-l-1<J:=\f{1}{2}\log_2|\bOm|,~~l=0,\cdots,L-1,~\aal\in\BB,~\bk\in\bOm,
\end{align*}
which shows that our choice of $\lam$ satisfies A2. When the wavelet frame transform is applied to $\bsf$ in \cref{Samples}, the refinable equations \cref{MRA-RF,MRA-Fra} and the definitions of $\msX^J(\Ps)$ \cref{QASystem,QAFramelet} lead to
\begin{align*}
\big(\bmW_{l,\aal}\bsf\big)[\bk]=2^J\la f,\psi_{\aal,J-1-l,\bk}\ra,~~~~~\aal\in\BB\cup\left\{\0\right\},~~~0\leq l\leq L-1
\end{align*}
with $\psi_{\0,l,\bk}=\phi_{l,\bk}$ for simplicity. Then together with \cref{DecayCondition}, we have
\begin{align*}
\left\|\lam\cdot\bmW\bsf\right\|_1\leq 2^J\sum_{l=J-L}^{\infty}2^{\beta l}\sum_{\bk\in\bOm}\sum_{\aal\in\BB}\big|\la f,\psi_{\aal,l,\bk}\ra\big|=C_f2^J,
\end{align*}
which shows that $\left\|\lam\cdot\bmW\bsf\right\|_1\leq C_{\bsf}2^J$ with $C_{\bsf}=C_f$. Hence, since we have $\left\|\bsf\right\|_{\ell_{\infty}(\bOm)}\leq\left\|f\right\|_{L_{\infty}(\Om)}\leq M$ by H\"older's inequality (e.g. \cite{Folland1999}), \cref{CorollaryResult1} follows from \cref{Th2} by setting
\begin{align*}
\wt{c}=\f{64M^2\left\|\bmA\right\|_{\infty}^2}{3\s_{\min}^2(\bmA)}\left[4+3\sqrt{5(4a+\max\{1-\beta,1\})C_{\bmW}C_f(1+C_{\bmW}C_f)}\right]\sqrt{2},
\end{align*}
for some $a\geq1$, which is independent of $J$, $L$, $\rho$, and $\eta$.

For \cref{CorollaryResult2}, note that we have
\begin{align*}
\left\|u_J^{\La}-f\right\|_{L_2(\Om)}\leq\left\|u_J^{\La}-f_J\right\|_{L_2(\Om)}+\left\|f_J-f\right\|_{L_2(\Om)}\leq\f{n}{2^J}\left\|\bu^{\La}-\bsf\right\|_{\ell_2(\bOm)}+\left\|f_J-f\right\|_{L_2(\Om)},
\end{align*}
where the last inequality comes from \cref{AdjointBessel} in \cref{Lemma3}, and $n$ denotes the order of $\phi$ (the tensor product B-spline used). Hence, we have
\begin{align*}
\left\|u_J^{\La}-f\right\|_{L_2(\Om)}^2\leq2\left(\f{n^2}{2^{2J}}\left\|\bu^{\La}-\bsf\right\|_{\ell_2(\bOm)}^2+\left\|f_J-f\right\|_{L_2(\Om)}^2\right),
\end{align*}
and by \cref{CorollaryResult1}, the proof is completed when we estimate $\left\|f_J-f\right\|_{L_2(\Om)}^2$. The standard wavelet frame decomposition gives
\begin{align*}
f=\sum_{\bk\in\bOm}\la f,\phi_{J,\bk}\ra\phi_{J,\bk}+\sum_{l=J}^{\infty}\sum_{\bk\in\bOm}\sum_{\aal\in\BB}\la f,\psi_{\aal,l,\bk}\ra\psi_{\aal,l,\bk}.
\end{align*}
Then the Bessel property of $\msX^J(\Ps)$ gives
\begin{align*}
\left\|f-f_J\right\|_{L_2(\Om)}^2\leq\sum_{l=J}^{\infty}\sum_{\bk\in\bOm}\sum_{\aal\in\BB}\big|\la f,\psi_{\aal,l,\bk}\ra\big|^2.
\end{align*}
Note that by H\"older's inequality, we have
\begin{align*}
\big|\la f,\psi_{\aal,l,\bk}\ra\big|\leq\left\|f\right\|_{L_{\infty}(\Om)}\left\|\psi_{\aal,l,\bk}\right\|_{L_1(\R^2)}=2^{-l}\left\|f\right\|_{L_{\infty}(\Om)}\left\|\psi_{\aal}\right\|_{L_1(\R^2)}\leq 2^{-l}M\left\|\psi_{\aal}\right\|_{L_1(\R^2)},
\end{align*}
for all $l\geq J$, and this leads to
\begin{align*}
\left\|f-f_J\right\|_{L_2(\Om)}^2&\leq M\left(\max_{\aal\in\BB}\left\|\psi_{\aal}\right\|_{L_1(\R^2)}\right)\sum_{l=J}^{\infty}2^{-l}\sum_{\bk\in\bOm}\sum_{\aal\in\BB}\big|\la f,\psi_{\aal,l,\bk}\ra\big|\\
&\leq M\left(\max_{\aal\in\BB}\left\|\psi_{\aal}\right\|_{L_1(\R^2)}\right)\sum_{l=J}^{\infty}2^{-l}\f{2^{(\beta+1)l}}{2^{(\beta+1)J}}\sum_{\bk\in\bOm}\sum_{\aal\in\BB}\big|\la f,\psi_{\aal,l,\bk}\ra\big|\\
&\leq M\left(\max_{\aal\in\BB}\left\|\psi_{\aal}\right\|_{L_1(\R^2)}\right)2^{-(\beta+1)J}\sum_{l=J-L}^{\infty}2^{\beta l}\sum_{\bk\in\bOm}\sum_{\aal\in\BB}\big|\la f,\psi_{\aal,l,\bk}\ra\big|\\
&=C_fM\left(\max_{\aal\in\BB}\left\|\psi_{\aal}\right\|_{L_1(\R^2)}\right)2^{-(\beta+1)J}.
\end{align*}
Therefore, by letting
\begin{align*}
C_1=2n^2\wt{c},~~C_2=\f{32}{3}\s_{\min}^{-2}(\bmA)n^2,~~\text{and}~~C_3=2C_fM\left(\max_{\aal\in\BB}\left\|\psi_{\aal}\right\|_{L_1(\R^2)}\right)
\end{align*}
all of which are independent of $J$, $L$, $\rho$, and $\eta$, we obtain \cref{CorollaryResult2}, and this concludes \cref{Th3}.\qquad$\square$
\end{poth3}

\section{Technical proofs}\label{TechnicalProofs}

This section is devoted to the technical details left in \cref{ErrorAnalysis}. More precisely, we present the proofs of \cref{Th1}, \cref{Lemma1}, and \cref{Th4}.

\subsection{Proof of \cref{Th1}}\label{ProofTh1}

To prove \cref{Th1}, we need some technical lemmas and propositions. The idea follows the same line as in \cite[Theorem 2.3.]{J.F.Cai2011}, which is mainly based on the idea of statistical learning \cite{F.Cucker2002}. However, since our problem involves a linear operator $\bmA\neq\bmI$, we still include the proof for the completeness. We first introduce the following ratio probability inequality, which can be derived from Bernstein inequality \cite{F.Cucker2002,Vapnik1998}.

\begin{lemma}\label{Lemma2} Assume that a random variable $\xi$ on a probability space $\msX$ satisfies $\mu:=\EE\xi\geq0$ and $\left|\xi-\mu\right|\leq B$ almost surely. If $\EE\xi^2\leq c\EE\xi$ for some $c>0$, the following inequality
\begin{align*}
\PP\left\{\z\in\msX^m:\f{\mu-\f{1}{m}\sum_{i=1}^m\xi(z_i)}{\sqrt{\mu+\ep}}>\gamma\sqrt{\ep}\right\}\leq\exp\left\{-\f{3\gamma^2m\ep}{6c+2B}\right\}
\end{align*}
holds for arbitrary $\ep>0$ and $0<\gamma\leq1$.
\end{lemma}

Let $\zze$ be a random variable which is drawn from the uniform distribution on $\bOm$. For $\bu\in\msM$, we introduce $\xi=\left|(\bmA\bu)[\zze]-(\bmA\bsf)[\zze]\right|^2$. Then $\xi$ is a random variable which is drawn from the uniform distribution on $\bOm$, and satisfies $0\leq\xi\leq\left\|\bmA\right\|_{\infty}^2M^2$. Define
\begin{align} \label{FullExpect}
\mE(\bu)=\EE\xi=\EE\left(\left|(\bmA\bu)[\zze]-(\bmA\bsf)[\zze]\right|^2\right)=\f{1}{|\bOm|}\sum_{\bk\in\bOm}\left|(\bmA\bu)[\bk]-(\bmA\bsf)[\bk]\right|^2.
\end{align}
We also define $\mE_{\La}(\bu)$ as a conditional expectation of $\xi$ given $\La$:
\begin{align}\label{PartialExpect}
\mE_{\La}(\bu)=\EE\left(\xi\big|\La\right)=\f{1}{|\La|}\sum_{\bk\in\La}\left|(\bmA\bu)[\bk]-(\bmA\bsf)[\bk]\right|^2.
\end{align}
Then we have the following lemma.
\begin{lemma}\label{Lemma4} Let $\msM$ be defined as \cref{HypothesisSpace}. For $\bu_1,\bu_2\in\msM$, we have
\begin{align}
\left|\mE(\bu_1)-\mE(\bu_2)\right|&\leq2\left\|\bmA\right\|_{\infty}^2M\left\|\bu_1-\bu_2\right\|_{\ell_{\infty}(\bOm)}\label{Ineq:E}\\
\left|\mE_{\La}(\bu_1)-\mE_{\La}(\bu_2)\right|&\leq2\left\|\bmA\right\|_{\infty}^2M\left\|\bu_1-\bu_2\right\|_{\ell_{\infty}(\bOm)}\label{Ineq:E_La}.
\end{align}
\end{lemma}

\begin{pf} We have
\begin{align*}
\left|\mE(\bu_1)-\mE(\bu_2)\right|&=\f{1}{|\bOm|}\left|\sum_{\bk\in\bOm}\left|(\bmA\bu_1)[\bk]-(\bmA\bsf)[\bk]\right|^2-\sum_{\bk\in\bOm}\left|(\bmA\bu_2)[\bk]-(\bmA\bsf)[\bk]\right|^2\right|\\
&\leq\f{1}{|\bOm|}\sum_{\bk\in\bOm}\left|(\bmA\bu_1)[\bk]+(\bmA\bu_2)[\bk]-2(\bmA\bsf)[\bk]\right|\left|(\bmA\bu_1)[\bk]-(\bmA\bu_2)[\bk]\right|\\
&\leq\f{1}{|\bOm|}\left\|\bmA\left(\bu_1+\bu_2-2\bsf\right)\right\|_{\ell_{\infty}(\bOm)}\left\|\bmA\left(\bu_1-\bu_2\right)\right\|_{\ell_1(\bOm)}\\
&\leq\left\|\bmA\left(\bu_1+\bu_2-2\bsf\right)\right\|_{\ell_{\infty}(\bOm)}\left\|\bmA\left(\bu_1-\bu_2\right)\right\|_{\ell_{\infty}(\bOm)}\\
&\leq2\left\|\bmA\right\|_{\infty}^2M\left\|\bu_1-\bu_2\right\|_{\ell_{\infty}(\bOm)},
\end{align*}
where we have used the fact that $\bu_1,\bu_2,\bsf\in[0,M]^{\bOm}$ in the last inequality; that is, $\left\|\bu_1+\bu_2-2\bsf\right\|_{\ell_{\infty}(\bOm)}\leq2M$. \cref{Ineq:E_La} can be proved in the similar way by replacing $\bOm$ with $\La$. \qquad$\square$
\end{pf}

With the aid of \cref{Lemma2,Lemma4}, we can present the following ratio probability inequality involving the space $\msM$ defined as \cref{HypothesisSpace}.

\begin{proposition}\label{Prop2} Let $\msM$, $\mE(\bu)$ and $\mE_{\La}(\bu)$ be respectively defined as \cref{HypothesisSpace}, \cref{FullExpect}, and \cref{PartialExpect}. Then for any $\ep>0$ and $0<\gamma<1$, we have
\begin{align*}
\PP\left\{\sup_{\bu\in\msM}\f{\mE(\bu)-\mE_{\La}(\bu)}{\sqrt{\mE(\bu)+\ep}}>4\gamma\sqrt{\ep}\right\}\leq\mN\left(\msM,\f{\gamma\ep}{2\left\|\bmA\right\|_{\infty}^2M}\right)\exp\left\{-\f{3\gamma^2m\ep}{8\left\|\bmA\right\|_{\infty}^2M^2}\right\}
\end{align*}
where $\mN(\msM,r)$ is the covering number of $\msM$ with respect to $\left\|\cdot\right\|_{\ell_{\infty}(\bOm)}$.
\end{proposition}

\begin{pf} Let $K=\mN\left(\msM,\f{\gamma\ep}{2\left\|\bmA\right\|_{\infty}^2M}\right)$ and we choose $\left\{\bu_1,\cdots,\bu_K\right\}\subseteq\msM$ such that
\begin{align*}
\msM\subseteq\bigcup_{j=1}^K\left\{\bu:\left\|\bu-\bu_j\right\|_{\ell_{\infty}(\bOm)}\leq\f{\gamma\ep}{2\left\|\bmA\right\|_{\infty}^2M}\right\}.
\end{align*}
For each $j=1,\cdots,K$, we consider the random variable $\xi=\left|(\bmA\bu_j)[\zze]-(\bmA\bsf)[\zze]\right|^2$ where $\zze$ is a random variable drawn from the uniform distribution on $\bOm$. Note that $\xi$ is a random variable on $\bOm$, and
\begin{align*}
\mu=\EE\xi=\f{1}{|\bOm|}\sum_{\bk\in\bOm}\left|(\bmA\bu_j)[\bk]-(\bmA\bsf)[\bk]\right|^2=\mE(\bu_j).
\end{align*}
Then since $\bu_j$, $\bsf\in\msM$, we have $0\leq\xi\leq\left\|\bmA\right\|_{\infty}^2M^2$, $\left|\xi-\EE\xi\right|\leq\left\|\bmA\right\|_{\infty}^2M^2$, and
\begin{align*}
\EE\xi^2=\EE\left(\left|(\bmA\bu_j)[\zze]-(\bmA\bsf)[\zze]\right|^4\right)\leq\left\|\bmA\right\|_{\infty}^2M^2\EE\left(\left|(\bmA\bu_j)[\zze]-(\bmA\bsf)[\zze]\right|^2\right)=\left\|\bmA\right\|_{\infty}^2M^2\EE\xi.
\end{align*}
In addition,
\begin{align*}
\f{1}{m}\sum_{\bk\in\La}\xi(\bk)=\f{1}{|\La|}\sum_{\bk\in\La}\left|(\bmA\bu_j)[\bk]-(\bmA\bsf)[\bk]\right|^2=\mE_{\La}(\bu_j).
\end{align*}
Applying \cref{Lemma2} to $\xi$ with $B=c=\left\|\bmA\right\|_{\infty}^2M^2$, we have
\begin{align}\label{RecallLemma2}
\PP\left\{\f{\mE(\bu_j)-\mE_{\La}(\bu_j)}{\sqrt{\mE(\bu_j)+\ep}}>\gamma\sqrt{\ep}\right\}\leq\exp\left\{-\f{3\gamma^2m\ep}{8\left\|\bmA\right\|_{\infty}^2M^2}\right\}.
\end{align}
For each $\bu\in\msM$, we can find $\bu_j$ such that $\left\|\bu-\bu_j\right\|_{\infty}\leq\f{\gamma\ep}{2\left\|\bmA\right\|_{\infty}^2M}$. Then by \cref{Lemma4}, we have
\begin{align*}
\left|\mE_{\La}(\bu)-\mE_{\La}(\bu_j)\right|&\leq2\left\|\bmA\right\|_{\infty}^2M\left\|\bu-\bu_j\right\|_{\ell_{\infty}(\bOm)}\leq\gamma\ep,\\
\left|\mE(\bu)-\mE(\bu_j)\right|&\leq2\left\|\bmA\right\|_{\infty}^2M\left\|\bu-\bu_j\right\|_{\ell_{\infty}(\bOm)}\leq\gamma\ep.
\end{align*}
In other words,
\begin{align*}
\f{\left|\mE_{\La}(\bu)-\mE_{\La}(\bu_j)\right|}{\sqrt{\mE(\bu)+\ep}}\leq\gamma\sqrt{\ep}~~\text{and}~~\f{\left|\mE(\bu)-\mE(\bu_j)\right|}{\sqrt{\mE(\bu)+\ep}}\leq\gamma\sqrt{\ep}.
\end{align*}
In particular, the second inequality implies that
\begin{align*}
\mE(\bu_j)+\ep&=\mE(\bu_j)-\mE(\bu)+\mE(\bu)+\ep\leq\gamma\sqrt{\ep}\sqrt{\mE(\bu)+\ep}+\mE(\bu)+\ep\\
&\leq\sqrt{\ep}\sqrt{\mE(\bu)+\ep}\leq\f{\mE(\bu)+2\ep}{2}+\mE(\bu)+\ep\leq2\left(\mE(\bu)+\ep\right).
\end{align*}
Hence, for an arbitrary $\bu\in\msM$ with $\left\|\bu-\bu_j\right\|_{\ell_{\infty}(\bOm)}\leq\f{\gamma\ep}{2\left\|\bmA\right\|_{\infty}^2M}$, we have $\sqrt{\mE(\bu_j)+\ep}\leq2\sqrt{\mE(\bu)+\ep}$.

Given that $\f{\mE(\bu)-\mE_{\La}(\bu)}{\sqrt{\mE(\bu)+\ep}}>4\gamma\sqrt{\ep}$, we have
\begin{align*}
\f{\mE(\bu_j)-\mE_{\La}(\bu_j)}{2\sqrt{\mE(\bu)+\ep}}&\geq\f{\mE(\bu_j)-\mE(\bu)}{2\sqrt{\mE(\bu)+\ep}}+\f{\mE(\bu)-\mE_{\La}(\bu)}{2\sqrt{\mE(\bu)+\ep}}+\f{\mE_{\La}(\bu)-\mE_{\La}(\bu_j)}{2\sqrt{\mE(\bu)+\ep}}\\
&>-\f{\gamma\sqrt{\ep}}{2}+2\gamma\sqrt{\ep}-\f{\gamma\sqrt{\ep}}{2}=\gamma\sqrt{\ep}.
\end{align*}
Then since we have $\sqrt{\mE(\bu_j)+\ep}\leq2\sqrt{\mE(\bu)+\ep}$ for any $\bu\in\msM$ with $\left\|\bu-\bu_j\right\|_{\infty}\leq\f{\gamma\ep}{2\left\|\bmA\right\|_{\infty}^2M}$, it follows that if $\f{\mE(\bu)-\mE_{\La}(\bu)}{\sqrt{\mE(\bu)+\ep}}>4\gamma\sqrt{\ep}$ holds, then the following inequality
\begin{align*}
\f{\mE(\bu_j)-\mE_{\La}(\bu_j)}{\sqrt{\mE(\bu_j)+\ep}}\geq\f{\mE(\bu_j)-\mE_{\La}(\bu)}{2\sqrt{\mE(\bu)+\ep}}>\gamma\sqrt{\ep}
\end{align*}
holds. Hence, for each fixed $j=1,\cdots,K$,
\begin{align*}
\PP\left\{\sup_{\bu\in\{\left\|\bu-\bu_j\right\|_{\ell_{\infty}(\bOm)}\leq\f{\gamma\ep}{2\left\|\bmA\right\|_{\infty}^2M}\}}\f{\mE(\bu)-\mE_{\La}(\bu)}{\sqrt{\mE(\bu)+\ep}}>4\gamma\sqrt{\ep}\right\}\leq\PP\left\{\f{\mE(\bu_j)-\mE_{\La}(\bu_j)}{\sqrt{\mE(\bu_j)+\ep}}>\gamma\sqrt{\ep}\right\}.
\end{align*}
Then since $\msM\subseteq\cup_{j=1}^K\left\{\bu:\left\|\bu-\bu_j\right\|_{\ell_{\infty}(\bOm)}\leq\f{\gamma\ep}{2\left\|\bmA\right\|_{\infty}^2M}\right\}$, together with \cref{RecallLemma2}, we have
\begin{align*}
\PP\left\{\sup_{\bu\in\msM}\f{\mE(\bu)-\mE_{\La}(\bu)}{\sqrt{\mE(\bu)+\ep}}>4\gamma\sqrt{\ep}\right\}&\leq\sum_{j=1}^K\PP\left\{\sup_{\bu\in\{\left\|\bu-\bu_j\right\|_{\ell_{\infty}(\bOm)}\leq\f{\gamma\ep}{2\left\|\bmA\right\|_{\infty}^2M}\}}\f{\mE(\bu)-\mE_{\La}(\bu)}{\sqrt{\mE(\bu)+\ep}}>4\gamma\sqrt{\ep}\right\}\\
&\leq\sum_{j=1}^K\PP\left\{\f{\mE(\bu_j)-\mE_{\La}(\bu_j)}{\sqrt{\mE(\bu_j)+\ep}}>\gamma\sqrt{\ep}\right\}\\
&\leq\mN\left(\msM,\f{\gamma\ep}{2\left\|\bmA\right\|_{\infty}^2M}\right)\exp\left\{-\f{3\gamma^2m\ep}{8\left\|\bmA\right\|_{\infty}^2M^2}\right\}.
\end{align*}
This concludes \cref{Prop2}.\qquad$\square$
\end{pf}

Now, we can provide the proof of \cref{Th1}.

\begin{poth1} By \cref{Prop2}, for arbitrary $\ep>0$ and $0<\gamma<1$, the inequality
\begin{align*}
\sup_{\bu\in\msM}\f{\mE(\bu)-\mE_{\La}(\bu)}{\sqrt{\mE(\bu)+\ep}}\leq4\gamma\sqrt{\ep}
\end{align*}
holds with probability at least
\begin{align*}
1-\mN\left(\msM,\f{\gamma\ep}{2\left\|\bmA\right\|_{\infty}^2M}\right)\exp\left\{-\f{3\gamma^2m\ep}{8\left\|\bmA\right\|_{\infty}^2M^2}\right\}.
\end{align*}
Therefore, for all $\bu\in\msM$, the inequality
\begin{align}\label{Ineq:E_Modified}
\mE(\bu)-\mE_{\La}(\bu)\leq4\gamma\sqrt{\ep}\sqrt{\mE(\bu)+\ep}
\end{align}
holds with the same probability. Then by taking $\gamma=\sqrt{2}/8$ and $\bu=\bu^{\La}\in\msM$ in \cref{Ineq:E_Modified}, it follows that
\begin{align}\label{Ineq:Eu^La}
\mE\left(\bu^{\La}\right)-\mE_{\La}\left(\bu^{\La}\right)\leq\f{1}{2}\sqrt{2\ep\left(\mE\left(\bu^{\La}\right)+\ep\right)}
\end{align}
holds with probability at least $1-\mN\left(\msM,\f{\ep}{8\sqrt{2}\left\|\bmA\right\|_{\infty}^2M}\right)\exp\left\{-\f{3m\ep}{256\left\|\bmA\right\|_{\infty}^2M^2}\right\}$. Besides, since
\begin{align}\label{Ineq:E_Lau^La}
\mE_{\La}\left(\bu^{\La}\right)\leq\f{2}{|\La|}\sum_{\bk\in\La}\left|(\bmA\bu^{\La})[\bk]-\bsg[\bk]\right|^2+\f{2}{|\La|}\sum_{\bk\in\La}\left|\bsg[\bk]-(\bmA\bsf)[\bk]\right|^2\leq4\eta^2,
\end{align}
combining \cref{Ineq:Eu^La} and \cref{Ineq:E_Lau^La} yields
\begin{align*}
\msE\left(\bu^{\La}\right)\leq\f{1}{2}\sqrt{2\ep\left(\mE\left(\bu^{\La}\right)+\ep\right)}+4\eta^2.
\end{align*}
Noting that $\mN\left(\msM,\f{\ep}{8\sqrt{2}\left\|\bmA\right\|_{\infty}^2M}\right)\leq\mN\left(\msM,\f{\ep}{12\left\|\bmA\right\|_{\infty}^2M}\right)$, we have
\begin{align*}
\PP\left\{\mE\left(\bu^{\La}\right)\leq\ep+\f{16}{3}\eta^2\right\}\geq1-\mN\left(\msM,\f{\ep}{12\left\|\bmA\right\|_{\infty}^2M}\right)\exp\left\{-\f{3m\ep}{256\left\|\bmA\right\|_{\infty}^2M^2}\right\}.
\end{align*}
Finally, since
\begin{align*}
\f{1}{|\bOm|}\left\|\bu^{\La}-\bsf\right\|_{\ell_2(\bOm)}^2\leq\f{\s_{\min}^{-2}(\bmA)}{|\bOm|}\sum_{\bk\in\bOm}\left|(\bmA\bu^{\La})[\bk]-(\bmA\bsf)[\bk]\right|^2=\s_{\min}^{-2}(\bmA)\mE\left(\bu^{\La}\right),
\end{align*}
we therefore have
\begin{align*}
\PP\left\{\f{1}{|\bOm|}\left\|\bu^{\La}-\bsf\right\|_{\ell_2(\bOm)}^2\leq\ep+\f{16}{3}\s_{\min}^{-2}(\bmA)\eta^2\right\}&\geq\PP\left\{\mE\left(\bu^{\La}\right)\leq\s_{\min}^2(\bmA)\ep+\f{16}{3}\eta^2\right\}\\
&\geq1-\mN\left(\msM,\f{\s_{\min}^2(\bmA)\ep}{12\left\|\bmA\right\|_{\infty}^2M}\right)\exp\left\{-\f{3m\s_{\min}^2(\bmA)\ep}{256\left\|\bmA\right\|_{\infty}^2M^2}\right\}.
\end{align*}
This concludes \cref{Th1}. \qquad$\square$
\end{poth1}

\subsection{Proof of \cref{Lemma1}}\label{ProofLemma1}

Before we prove \cref{Lemma1}, we first mention that the authors in \cite{J.F.Cai2011} have presented the similar results for the generic MRA based wavelet frame systems, using the convergence rate of a stationary subdivision algorithm \cite{A.S.Cavaretta1991,Ji1995,W.Lawton1998}. However, since our definition of $\ell_1$-norm of wavelet frame coefficients does not include the low-pass filter coefficients, our proof is based on the technique different from \cite{J.F.Cai2011}. In fact, since we restrict ourselves to the tensor product B-spline wavelet frame systems, our proof follows the line described in \cite{J.Yang2017}, based on the definition of B-spline wavelet frame system in \cite{A.Ron1997}.

We denote respectively by $\bp=\left[\begin{array}{cc}
1&1
\end{array}\right]$ one dimensional $1$st order average filter and by $\bq=\left[\begin{array}{cc}
1&-1
\end{array}\right]$ one dimensional $1$st order difference filter. Then we note that for $\bu\in\R^{\bOm}$, the direct computation gives
\begin{align}\label{TVUB1}
\left\|\na\bu\right\|_1\leq\left\|\bq_{\e_1}[-\cdot]\circledast\bu\right\|_{\ell_1(\bOm)}+\left\|\bq_{\e_2}[-\cdot]\circledast\bu\right\|_{\ell_1(\bOm)}
\end{align}
where $\left\|\Na\bu\right\|_1$ is defined as \cref{DiscreteTV}, and the convolution $\circledast$ is defined as \cref{ConvolutionMeaning}. In \cref{TVUB1}, $\bq_{\e_1}$ and $\bq_{\e_2}$ with the standard basis $\left\{\e_1,\e_2\right\}$ for $\R^2$ is defined as
\begin{align}\label{DiscreteGradRev}
\bq_{\e_1}[\bk]=\bq[k_1]\dde[k_2]~~~\text{and}~~~\bq_{\e_2}[\bk]=\dde[k_1]\bq[k_2]
\end{align}
where $\dde$ is the one dimensional discrete Dirac delta; $\dde[k]=1$ for $k=0$, and $\dde[k]=0$ for $k\neq0$.

Now, we note that for $L\geq 1$, we have
\begin{align*}
\sum_{\bk\in\bOm}\sum_{\aal\in\BB}\left|\big(\bmW_{0,\aal}\bu\big)[\bk]\right|\leq\sum_{\bk\in\bOm}\sum_{l=0}^{L-1}\sum_{\aal\in\BB}\left|\big(\bmW_{l,\aal}\bu\big)[\bk]\right|.
\end{align*}
Hence, \cref{TVUB1} tells us that, if we find $C>0$ independent of $J$ such that
\begin{align}\label{Lemma1:Goal1}
\left\|\bq_{\e_1}[-\cdot]\circledast\bu\right\|_{\ell_1(\bOm)}+\left\|\bq_{\e_2}[-\cdot]\circledast\bu\right\|_{\ell_1(\bOm)}\leq C\sum_{\bk\in\bOm}\sum_{\aal\in\BB}\left|\big(\bmW_{0,\aal}\bu\big)[\bk]\right|,
\end{align}
we obtain \cref{Lemma1Result}, by simply setting $C_{\bmW}=\max\left\{C,1\right\}$.

Recall that the first equality of UEP \cref{UEP} is rewritten as
\begin{align*}
\sum_{\aal\in\BB\cup\{\0\}}\left(\a_{\aal}\ast\a_{\aal}[-\cdot]\right)[\bk]=\dde[\bk]=\left\{\begin{array}{rl}
1~&\text{if}~\bk=\0;\\
0~&\text{if}~\bk\neq\0,
\end{array}\right.
\end{align*}
where $\ast$ is the standard discrete convolution on $\Z^2$. Then for $j=1,2$, we have
\begin{align*}
\bq_{\e_j}[-\cdot]\circledast\bu=\sum_{\aal\in\BB\cup\{\0\}}\left(\a_{\aal}\ast\a_{\aal}[-\cdot]\right)\ast\bq_{\e_j}[-\cdot]\circledast\bu,
\end{align*}
and by the discrete Young's inequality (e.g. \cite{Beckner1975}), we have
\begin{align}\label{IntermediateResult}
\left\|\bq_{\e_j}[-\cdot]\circledast\bu\right\|_{\ell_1(\bOm)}\leq\sum_{\aal\in\BB\cup\{\0\}}\left\|\a_{\aal}\right\|_{\ell_1(\Z^2)}\left\|\a_{\aal}[-\cdot]\ast\bq_{\e_j}[-\cdot]\circledast\bu\right\|_{\ell_1(\bOm)}.
\end{align}
From the definition of B-spline wavelet frame system, the Fourier series of one dimensional filters $\a_{\alpha}$ is
\begin{align}
\wh{\a}_{\0}(\xi)&=e^{-ij_r\f{\xi}{2}}\cos^r\left(\f{\xi_1}{2}\right),\label{FourierFilterLow}\\
\wh{\a}_{\alpha}(\xi)&=-i^{\alpha}e^{-ij_r\f{\xi}{2}}\sqrt{\binom{r}{\alpha}}\cos^{r-\alpha}\left(\f{\xi}{2}\right)\sin^{\alpha}\left(\f{\xi}{2}\right),~~~~~\alpha=1,\cdots,r\label{FourierFilterHigh}
\end{align}
where $j_r=1$ if $r$ is odd; $j_r=0$ if $r$ is even. Noting that the Fourier series of $\bp$ and $\bq$ are
\begin{align*}
\wh{\bp}(\xi)=2e^{-i\f{\xi}{2}}\cos\left(\f{\xi}{2}\right),~~~~~\text{and}~~~~~\wh{\bq}(\xi)=-2ie^{i\f{\xi}{2}}\sin\left(\f{\xi}{2}\right),
\end{align*}
we define $\bp^{\alpha}$ and $\bq^{\alpha}$ with $\alpha\in\N$ as
\begin{align*}
\bp^{\alpha}=\underbrace{\bp\ast\bp\ast\cdots\ast\bp}_{\alpha\text{ times}}~~~~\text{and}~~~~\bq^{\alpha}=\underbrace{\bq\ast\bq\ast\cdots\ast\bq}_{\alpha\text{ times}},
\end{align*}
and we set $\bp^0=\bq^0=\dde$, the one dimensional discrete Dirac delta. Then by \cref{FourierFilterLow,FourierFilterHigh}, we have
\begin{align}
\a_0[k]&=2^{-r}\bp^r[k-j_r] \label{1DFilterExplicitLow}\\
\a_{\alpha}[k]&=(-1)^{\alpha+1}2^{-r}\sqrt{\binom{r}{\alpha}}\bp^{r-\alpha}\ast\bq^{\alpha}[k-j_r],~~\alpha=1,\cdots,r.\label{1DFilterExplicitHigh}
\end{align}
Using \cref{1DFilterExplicitLow,1DFilterExplicitHigh} and the discrete Young's inequality, we have
\begin{align*}
\left\|\a_{\0}\right\|_{\ell_1(\Z^2)}=1~~~~~\text{and}~~~~~\left\|\a_{\aal}\right\|_{\ell_1(\Z^2)}\leq\sqrt{\binom{\bsr}{\aal}}~~~~\aal\in\BB,
\end{align*}
where $\bsr=(r,r)$, $\binom{\aal}{\bbe}=\binom{\alpha_1}{\beta_1}\binom{\alpha_2}{\beta_2}$ for multi-indices $\aal=(\alpha_1,\alpha_2)$, $\bbe=(\beta_1,\beta_2)$, and the first equality comes from the fact that $\wh{\a}_{\0}(\0)=\sum\a_{\0}[\bk]=1$. From $\binom{r}{0}=1$, \cref{IntermediateResult} becomes
\begin{align*}
\left\|\bq_{\e_j}[-\cdot]\circledast\bu\right\|_{\ell_1(\bOm)}\leq\sum_{\aal\in\BB\cup\{\0\}}\sqrt{\binom{\bsr}{\aal}}\left\|\a_{\aal}[-\cdot]\ast\bq_{\e_j}[-\cdot]\circledast\bu\right\|_{\ell_1(\bOm)}.
\end{align*}
It is obvious that, for $\aal\in\BB$, we have
\begin{align*}
\left\|\a_{\aal}[-\cdot]\ast\bq_{\e_j}[-\cdot]\circledast\bu\right\|_{\ell_1(\bOm)}\leq\left\|\bq_{\e_j}[-\cdot]\right\|_{\ell_1(\Z)}\left\|\a_{\aal}[-\cdot]\circledast\bu\right\|_{\ell_1(\bOm)}=2\left\|\bmW_{0,\aal}\bu\right\|_{\ell_1(\bOm)}
\end{align*}
by the discrete Young's inequality. Hence, we need to bound $\left\|\a_{\0}[-\cdot]\ast\bq_{\e_j}[-\cdot]\circledast\bu\right\|_{\ell_1(\bOm)}$ by using $\left\|\bmW_{0,\aal}\bu\right\|_{\ell_1(\bOm)}$ for some suitable $\aal\in\BB$. To do this, we first note that
\begin{align*}
\a_0\ast\bq[k]=r^{-\f{1}{2}}\bp\ast\a_1[k-j_r].
\end{align*}
Then by the discrete Young's inequality, we have
\begin{align*}
\left\|\a_{\0}[-\cdot]\ast\bq_{\e_j}[-\cdot]\circledast\bu\right\|_{\ell_1(\bOm)}\leq\f{2}{\sqrt{r}}\left\|\a_{\e_j}[-\cdot]\circledast\bu\right\|_{\ell_1(\bOm)}=\f{2}{\sqrt{r}}\left\|\bmW_{0,\e_j}\bu\right\|_{\ell_1(\bOm)}.
\end{align*}
Combining these altogether, we have
\begin{align*}
\left\|\bq_{\e_1}[-\cdot]\circledast\bu\right\|_{\ell_1(\bOm)}&+\left\|\bq_{\e_2}[-\cdot]\circledast\bu\right\|_{\ell_1(\bOm)}\\
&\leq4\sum_{\aal\in\BB}\sqrt{\binom{\bsr}{\aal}}\left\|\bmW_{0,\aal}\bu\right\|_{\ell_1(\bOm)}+\f{2}{\sqrt{r}}\left\|\bmW_{0,\e_1}\bu\right\|_{\ell_1(\bOm)}+\f{2}{\sqrt{r}}\left\|\bmW_{0,\e_2}\bu\right\|_{\ell_1(\bOm)}\\
&\leq6\sum_{\aal\in\BB}\sqrt{\binom{\bsr}{\aal}}\left\|\bmW_{0,\aal}\bu\right\|_{\ell_1(\bOm)},
\end{align*}
where the last inequality follows from the fact that $r\geq1$ and $\binom{r}{\alpha}\geq1$. Therefore, letting
\begin{align*}
C=6\max_{\aal\in\BB}\sqrt{\binom{\bsr}{\aal}},
\end{align*}
we obtain \cref{Lemma1:Goal1}, and this concludes \cref{Lemma1}.

\subsection{Proof of \cref{Th4}}\label{ProofTh4}

\cref{Th4} is to estimate the covering number. The proof follows the line similar to \cite[Theorem 2.4]{J.F.Cai2011}. However, since we improve \cite[Theorem 2.4]{J.F.Cai2011} by relaxing the constraint of the radius $r$, we include the detailed proof for the sake of completeness. It is obvious that if $\bu\in[0,M]^{\bOm}$, then $\left\|\bu\right\|_{\ell_{\infty}(\bOm)}\leq M$. In addition, if A2 is satisfied, together with the definition of $\msM$ in \cref{HypothesisSpace} and \cref{Lemma1}, we have
\begin{align*}
\left\|\Na\bu\right\|_1\leq C_{\bmW}\left\|\bmW\bu\right\|_1\leq C_{\bmW}2^{\max\{-\beta,0\}\max\{\Upsilon(l,\aal,\bk):(l,\aal,\bk)\in\Gamma\}}\left\|\lam\cdot\bmW\bu\right\|_1.
\end{align*}
Since $\Upsilon:\Gamma\to\N$ satisfies \cref{lambdaCondition}, for $\bu\in\msM$, we have
\begin{align*}
\left\|\Na\bu\right\|_1\leq C_{\bmW}|\bOm|^{\f{\max\{-\beta,0\}}{2}}\left\|\lam\cdot\bmW\bsf\right\|_1.
\end{align*}
In other words, $\msM\subseteq\wt{\msM}$ and $\mN\left(\msM,r\right)\leq\mN\left(\wt{\msM},r\right)$ where $\wt{\msM}$ is defined as \cref{RelaxedFeasibleSet}. Hence, it suffices to bound the covering number $\mN\left(\wt{\msM},r\right)$. In addition, it is easy to see that if there exists a finite set $F\subseteq\wt{\msM}$ such that
\begin{align*}
\wt{\msM}\subseteq\bigcup_{\bq\in F}\left\{\bu:\left\|\bu-\bq\right\|_{\ell_{\infty}(\bOm)}\leq r\right\},
\end{align*}
we have $\mN\left(\wt{\msM},r\right)\leq\left|F\right|$. What we need now is to construct an appropriate set $F$ by exploiting the specific structure of $\wt{\msM}$, so that $\left|F\right|$ has an appropriate upper bound.

For this purpose, let $\kappa=\lceil2M/r\rceil$, and we define
\begin{align}\label{QuantizedRange}
\msR=\left\{-\kappa r/2,(-\kappa+1)r/2,\cdots,\kappa r/2\right\}.
\end{align}
By \cite[Lemma 4.4]{J.F.Cai2011}, for each $\bu\in\wt{\msM}$, there exists $Q(\bu)\in\msR^{\bOm}$ such that
\begin{align*}
\left\|\bu-Q(\bu)\right\|_{\ell_{\infty}(\bOm)}\leq r/2~~~\text{and}~~~\left\|\Na\left(Q(\bu)\right)\right\|_1\leq\left\|\Na\bu\right\|_1.
\end{align*}
Let
\begin{align*}
\wt{F}=\left\{\bq\in\ell_{\infty}(\bOm):\bq=Q(\bu)~~\text{for some}~~\bu\in\wt{\msM}\right\}\subseteq\msR^{\bOm}.
\end{align*}
Notice that for each $\bq\in\wt{F}$, there may be more than one $\bu\in\wt{\msM}$ such that $\bq=Q(\bu)$.

For each $\bq\in\wt{F}$, choose $\bu_{\bq}\in\wt{\msM}$ such that $\left\|\bu_{\bq}-\bq\right\|_{\ell_{\infty}(\bOm)}\leq r/2$, and define $F=\left\{\bu_{\bq}:\bq\in\wt{F}\right\}$. For an arbitrary $\bu\in\wt{\msM}$, there exists $\bq\in\wt{F}$ such that $\left\|\bu-\bq\right\|_{\ell_{\infty}(\bOm)}\leq r/2$. This implies
\begin{align*}
\left\|\bu-\bu_{\bq}\right\|_{\ell_{\infty}(\bOm)}\leq\left\|\bu-\bq\right\|_{\ell_{\infty}(\bOm)}+\left\|\bq-\bu_{\bq}\right\|_{\ell_{\infty}(\bOm)}\leq r,
\end{align*}
by the definition of $\bu_{\bq}$. Therefore,
\begin{align*}
\wt{\msM}\subseteq\bigcup_{\bu_{\bq}\in F}\left\{\bu:\left\|\bu-\bu_{\bq}\right\|_{\ell_{\infty}(\bOm)}\leq r\right\}~~\text{and}~~\mN\left(\wt{\msM},r\right)\leq\left|F\right|\leq\left|\wt{F}\right|.
\end{align*}
Thus, the covering number $\mN\left(\wt{\msM},r\right)$ is bounded by any upper bound of $\left|\wt{F}\right|$. Notice that each $\bq\in\wt{F}$ is uniquely determined by $\Na\bq$ and $\bq[1,\cdots,1]$. Since $\wt{F}$ is a subset of $\msR^{\bOm}$, there are $2\kappa+1$ choices for $\bq[1,\cdots,1]$. It remains to count the number of choices in $\Na\bq$. Define
\begin{align*}
\Na\wt{F}=\left\{\Na\bq:\bq\in\wt{F}\right\}.
\end{align*}
Then we need to bound $\left|\Na\wt{F}\right|$.

To do this, we first consider the uniform upper bound of $\left\|\Na\bq\right\|_1$ for $\bq\in\wt{F}$. By the definition of $\wt{F}$ and \cite[Lemma 4.4]{J.F.Cai2011}, for each $\bq\in\wt{F}$, there exists $\bu\in\wt{\msM}$ such that $\left\|\Na\bq\right\|_1\leq\left\|\Na\bu\right\|_1$. Since $\left\|\lam\cdot\bmW\bsf\right\|_1\leq C_{\bsf}|\bOm|^{1/2}$ by assumption, we further have
\begin{align*}
\left\|\Na\bq\right\|_1\leq\left\|\Na\bu\right\|_1\leq C_{\bmW}C_{\bsf}\left|\bOm\right|^{b/2}
\end{align*}
with $b=\max\{1-\beta,1\}$ for notational simplicity. Hence, for all $\bq\in\wt{F}$, we have $\left\|\Na\bq\right\|_1\leq Kr/2$, where
\begin{align*}
K=\left\lceil\f{2C_{\bmW}C_{\bsf}\left|\bOm\right|^{b/2}}{r}\right\rceil.
\end{align*}
In addition, since $\bq\in\msR^{\bOm}$, each element of $\Na\bq$ has to be a multiple of $r/2$, which means that the range of $\Na\bq$ is a subset of $\left\{-Kr/2,-(K-1)r/2,\cdots,Kr/2\right\}^2$. Recall that there are $R=2\left(\left|\bOm\right|-\left|\bOm\right|^{1/2}\right)$ elements in $\Na\bq$. Hence, $\left|\Na\wt{F}\right|$ can be bounded by the number of all possible integer solutions of the following inequality
\begin{align*}
\left|x_1\right|+\left|x_2\right|+\cdots+\left|x_{R-1}\right|+\left|x_R\right|\leq K.
\end{align*}
That is (e.g. \cite{Ross1984}),
\begin{align*}
\left|\Na\wt{F}\right|&\leq1+\sum_{k=1}^K\sum_{s=1}^{\min\{k,R\}}2^s\binom{R}{R-s}\binom{k-1}{s-1}\\
&\leq1+\sum_{k=1}^K2^k\sum_{s=1}^{\min\{k,R\}}\binom{R}{R-s}\binom{k-1}{s-1}=1+\sum_{k=1}^K2^k\binom{R+k-1}{R-1}\\
&\leq\binom{R+K-1}{R-1}\sum_{k=0}^K2^k\leq2^{K+1}\binom{R+K-1}{K}\leq2\left[2\left(R+K-1\right)\right]^K.
\end{align*}
Hence, we have
\begin{align*}
\left|\wt{F}\right|\leq\left(4\kappa+2\right)\left[2\left(R+K-1\right)\right]^K.
\end{align*}
In other words, using $r\geq\left|\bOm\right|^{-a}$ with $a\geq1$ and $K-1\leq2C_{\bmW}C_{\bsf}|\bOm|^{a+b/2}$, we have
\begin{align*}
\ln\mN\left(\msM,r\right)&\leq\f{2C_{\bmW}C_{\bsf}\left|\bOm\right|^{b/2}}{r}\ln\left(2\left(R+K-1\right)\right)+\ln\left(\f{8M}{r}+2\right)\\
&\leq\f{2C_{\bmW}C_{\bsf}\left|\bOm\right|^{b/2}}{r}\left[\ln\left(4\left|\bOm\right|+\f{4C_{\bmW}C_{\bsf}\left|\bOm\right|^{b/2}}{r}\right)+\ln\f{10M}{r}\right]\\
&\leq\f{2C_{\bmW}C_{\bsf}\left|\bOm\right|^{b/2}}{r}\ln\left(\left(1+C_{\bmW}C_{\bsf}\right)40M\left|\bOm\right|^{2a+b/2}\right)
\end{align*}
where we use the fact that $a+b/2\geq1$ from the choice of $a$ and $b$ in the final inequality. Therefore, we have
\begin{align*}
\ln\mN\left(\msM,r\right)\leq\f{20M(4a+b)C_{\bmW}C_{\bsf}\left(1+C_{\bmW}C_{\bsf}\right)\left|\bOm\right|^{b/2}}{r}\log_2\left|\bOm\right|,
\end{align*}
and this completes the proof.\qquad$\square$

\section{Conclusion}\label{Conclusion}

In this paper, we present an approximation analysis of wavelet frame based restoration from degraded and incomplete measurements. By the combination of the uniform law of large numbers and the estimation for its involved covering number of a hypothesis space of the solution, we establish unified approximation properties of the wavelet frame based image restoration model, together with the improved estimate of covering number over the previous one in \cite{J.F.Cai2011}. In addition, thanks to the underlying multiresolution analysis structure, we further connect the error analysis in the discrete setting to the approximation of underlying function where the discrete data comes from. For the future works, we plan to study an approximation property of wavelet frame based image restoration from partial Fourier samples, where $\bmA=\bmsF$ (the unitary discrete Fourier transform) satisfies $\left\|\bmsF\right\|_{\infty}=|\bOm|^{1/2}$ and $\sigma_{\min}(\bmsF)=\sigma_{\max}(\bmsF)=1$ in our framework. We may also plan to establish an approximation of tight frame based missing data restoration on the graph \cite{Dong2017}, and the multiresolution approximation of functions on the manifold \cite{B.Dong2016}.

\section*{References}

\bibliographystyle{cas-model2-names}

\end{document}